\documentclass{article}%
\usepackage{amsmath}
\usepackage{graphicx}
\usepackage{epsfig}
\usepackage{multicol}
\usepackage{float}
\usepackage{amsfonts}
\usepackage{amssymb}
\usepackage{color}
\usepackage{xcolor}%
\newtheorem{theorem}{Theorem}[section]

\newtheorem{corollary}[theorem]{Corollary}

\newtheorem{lemma}[theorem]{Lemma}

\newtheorem{proposition}[theorem]{Proposition}
\newtheorem{remark}{Remark}[section]

\newenvironment{proof}[1][Proof]{\textbf{#1.} }{\ \rule{0.5em}{0.5em}}

\begin{document}

\title{Asymptotic behaviour of global solutions to a model of cell invasion}
\date{}
\author{{\normalsize Gabriela\ Li\c{t}canu* and Cristian Morales-Rodrigo**}
\and {\small *Institute of Mathematics "O. Mayer", Romanian
Academy,}\\{\small \ 700505 Ia\c{s}i, Romania. e-mail:
litcanu@uaic.ro}\\{\small **Dpto. de Ecuaciones Diferenciales y
An\'{a}lisis Num\'{e}rico, Universidad de Sevilla, } \\{\small Apdo.
de Correos 1160, 41080 Sevilla, Spain. e-mail: cristianm@us.es}}
\maketitle

\begin{abstract}
{\small In this paper we analyze a mathematical model focusing on key events
of the cells invasion process. Global well-possedness and asymptotic behaviour
of nonnegative solutions to the corresponding coupled system of three
nonlinear partial differential equations are studied.}

\


\ \

\textit{Key words or phrases}{\small : mathematical model, global existence, a
priori estimates, asymptotic behaviour}

\end{abstract}

\section{Introduction}

\label{introduction}

\setcounter{definition}{0}\setcounter{equation}{0}

In this paper we focus on a mathematical model describing the process of cells
invasion in the surrounding extracellular matrix. Because of the key role
played by the invasive processes in biological phenomena like, for example,
wound healing, morphogenesis or tumour invasion, there are a large number of
studies concerning them.

In \cite{byrne} the authors developed a mathematical model in order to
describe the migration of tumour cells through a collagen gel. More precisely,
the model is based on the hypothesis that the cells invasion is the final
result of the triad of adhesion, proteolysis and motility such that in contact
with the extracellular matrix, the invasive tumour cells produce proteolytic
enzymes which degrade it favoring the migration.

Our objective is to study a version of a model which underlies the models
proposed in \cite{chaplain} and \cite{byrne} (see also the references therein)
which involves three key variables: $u(x,t)$ the density of invasive cells,
$v(x,t)$ the density of extracellular matrix and $m(x,t)$ the concentration of
degradative enzymes such as proteases, each of them considered at $x\in\Omega$
and time $t>0$. Through this paper $\Omega\subset\mathbb{R}^{N}$,
$N\geqslant1$ is a bounded domain with a regular boundary. The model is the
following:%
\begin{align}
&  \frac{\partial u}{\partial t}=\underset{\text{\textit{diffusion}}%
}{\underbrace{d_{1}\Delta u}}-\underset{\text{\textit{haptotaxis}}%
}{\underbrace{\alpha_{1}\nabla\cdot(u\chi(v)\nabla v)}}+\underset
{\text{\textit{cells proliferation}}}{\underbrace{\alpha_{2}u\left(
1-\alpha_{3}u-\alpha_{4}v\right)  }} & x  &  \in\Omega,\quad t\in
\mathbb{R}_{+}\label{eq1_gen}\\
&  \frac{\partial v}{\partial t}=-\underset{\text{\textit{degradation}}%
}{\underbrace{\lambda mv}} & x  &  \in\Omega,\quad t\in\mathbb{R}%
_{+}\label{eq2_gen}\\
&  \frac{\partial m}{\partial t}=\underset{\text{\textit{diffusion}}%
}{\underbrace{d_{2}\Delta m}}-\underset{\text{\textit{decay}}}{\underbrace
{\beta_{1}m}}+\underset{\text{\textit{production}}}{\underbrace{\beta
_{2}ug(v)}} & x  &  \in\Omega,\quad t\in\mathbb{R}_{+} \label{eq3_gen}%
\end{align}
where the coefficients $d_{1}$, $d_{2}$, $\lambda$, $\alpha_{1}$, $\alpha_{3}%
$, $\alpha_{4}$, $\beta_{1}$, $\beta_{2}$ are positive constants and
$\alpha_{2}$ is non-negative. It is assumed that the change in time of the
cells density is due to the diffusion, to the haptotaxis with respect to
spatial gradients in the collagen gel and to the process involving
proliferation and degradation if $\alpha_{2}\neq0$. The new cells are created
at the rate $\alpha_{2}>0$ and their degradation is caused by the cells death
or neutralization due to the presence of the collagen gel. The function $\chi$
describes the sensitivity of the cells to spatial gradients of collagen and is
assumed to be non-negative.

In the equation (\ref{eq2_gen}) it is assumed that the proteases degrade the
collagen gel at the rate $\lambda$.

The proteases concentration is affected by the diffusion, natural decay and
production, this last term being proportional to the product of the cell
density and a non-negative function $g$ depending on the collagen gel concentration.

In what follows we consider the system (\ref{eq1_gen})-(\ref{eq3_gen})
together with the boundary conditions%
\begin{equation}
d_{1}\frac{\partial u}{\partial\eta}-\alpha_{1}u\chi(v)\frac{\partial
v}{\partial\eta}=\frac{\partial m}{\partial\eta}=0,\qquad x\in\partial
\Omega,\quad t\in\mathbb{R}_{+} \label{bc_gen}%
\end{equation}
where $\eta$ denotes the unit outward normal vector of $\partial\Omega$. We
have supposed that there is no flux of cells or proteases across the boundary
of the domain. The same hypothesis is assumed for the extracellular matrix.

We consider also the initial conditions%
\begin{equation}
(u,v,m)(x,0)=(u_{0},v_{0},m_{0})(x),\qquad x\in\Omega. \label{ic_gen}%
\end{equation}

Taking into account the biological interpretation for the solution of the
system (\ref{eq1_gen})-(\ref{eq3_gen}), we shall require that the functions
$u_{0}(x),$ $v_{0}(x),$ and $m_{0}(x)$ are non-negative.

For $g(v)=v$ and for a constant chemotactic coefficient $\chi(v)=\chi$, the
mathematical model (\ref{eq1_gen})-(\ref{ic_gen}) was previously studied
numerically in \cite{byrne} assuming that the cells invasion was radially
symmetrical. In the same hypotheses for $g$ and $\chi$, the system was
considered in \cite{tao1} where the existence of global solutions was
investigated in the 3-dimensional case. In this paper we prove the global
existence of the solutions when $g$ and $\chi$ are arbitrary functions
satisfying some hypotheses that will be given later. Moreover, we also study
the asymptotic behaviour of the solutions.

After nondimensionalizing the system (\ref{eq1_gen})-(\ref{eq3_gen}) and
redenoting the functions $\chi$ and $g$ if necessary, it becomes%
\begin{align}
&  \frac{\partial u}{\partial t}=\Delta u-\nabla\cdot(u\chi(v)\nabla v)+\mu
u\left(  1-u-v\right)  & x  &  \in\Omega,\quad t\in\mathbb{R}_{+}\label{eq1}\\
&  \frac{\partial v}{\partial t}=-mv & x  &  \in\Omega,\quad t\in
\mathbb{R}_{+}\label{eq2}\\
&  \frac{\partial m}{\partial t}=d\Delta m-\gamma m+ug(v) & x  &  \in
\Omega,\quad t\in\mathbb{R}_{+}\label{eq3}\\
&  \frac{\partial u}{\partial\eta}-u\chi(v)\frac{\partial v}{\partial\eta
}=\frac{\partial m}{\partial\eta}=0 & x  &  \in\partial\Omega,\quad
t\in\mathbb{R}_{+}\label{bc}\\
&  (u,v,m)(x,0)=(u_{0},v_{0},m_{0})(x) & x  &  \in\Omega, \label{ic}%
\end{align}
where%
\[
\mu=\left\{
\begin{array}
[c]{cc}%
1 & \text{, if \ }\alpha_{2}\neq0\\
0 & \text{, if \ }\alpha_{2}=0
\end{array}
\right.  ,\qquad d=\frac{d_{2}}{d_{1}},\qquad\gamma=\left\{
\begin{array}
[c]{cc}%
\beta_{1}\left(  \alpha_{2}\right)  ^{-1} & \text{, if \ }\alpha_{2}\neq0\\
1 & \text{, if \ }\alpha_{2}=0
\end{array}
\right.  .
\]
We notice that the new functions $\chi$ and $g$ are obtained from the initial
ones by a rescaling that does not change their initial properties.

Finally, we mention here that we can obtain an equivalent system to
(\ref{eq1})-(\ref{ic}) by making the change of variables%
\begin{equation}
w:=uz,\hspace{25pt}z:=e^{-\int_{0}^{v}\chi(s)ds}. \label{changevariables2}%
\end{equation}
In this way the system (\ref{eq1})-(\ref{ic}) transforms into%
\begin{align}
&  \frac{\partial w}{\partial t}=\Delta w+\chi(v)\nabla v\cdot\nabla w+\mu
w\left(  1-wz^{-1}-v\right)  +\chi(v)wvm & x  &  \in\Omega,\quad
t\in\mathbb{R}_{+}\label{e1}\\
&  \frac{\partial v}{\partial t}=-mv & x  &  \in\Omega,\quad t\in
\mathbb{R}_{+}\label{e2}\\
&  \frac{\partial m}{\partial t}=d\Delta m-\gamma m+wz^{-1}g(v) & x  &
\in\Omega,\quad t\in\mathbb{R}_{+}\label{e3}\\
&  \frac{\partial w}{\partial\eta}=\frac{\partial m}{\partial\eta}=0 & x  &
\in\partial\Omega,\quad t\in\mathbb{R}_{+}\label{b}\\
&  (w,v,m)(x,0)=(u_{0}e^{-\int_{0}^{v_{0}}\chi(s)ds},v_{0},m_{0}%
)(x)=(w_{0},v_{0},m_{0})(x) & x  &  \in\Omega. \label{i}%
\end{align}

In this paper we are concerning to prove the existence of a unique global
solution of (\ref{eq1})-(\ref{ic}) and also to investigate asymptotic
behaviour of the solution in a more general case when $d$ and $\gamma$ are
positive constants and $\mu$ is a non-negative constant.

We mention that in the case when $g(v)=v$ and $\chi(v)=\chi$ is a constant, in
\cite{tao1} the authors proved the global existence of solutions of
(\ref{e1})-(\ref{i}) in the 3-dimensional case based on a priori estimates. In
order to derive $L^{p}$ estimates, they used a similar approach as in
\cite{walker}.

We point out that in what follows we establish that the a-priori $L^{p}$
estimates, $p\geqslant1$, for the variable $w$ (and also for $u$) are uniform
in time. This fact will be important when we establish the asymptotic
behaviour of the solutions.

This paper is organized as follows. In Section \ref{preliminaries} we give the
notations and terminology used through the paper. In Section \ref{local} we
prove the local existence and the non-negativity of solutions for non-negative
initial data using a fixed point method. In Section \ref{global} we show that
the solution constructed in the previous section can be prolonged in time
until infinity when $N=3$. Section \ref{steadyst} is devoted to the stationary
problem associated to (\ref{e1})-(\ref{i}). In the last Section we show the
convergence to the steady-states and we obtain an explicit rate of convergence
in some cases.

\section{Preliminaries and notations}

\label{preliminaries}

\setcounter{definition}{0}\setcounter{equation}{0}

In this section we collect some tools and notations that will be used in the paper.

Let $\Omega\subset\mathbb{R}^{N}$, $N\geqslant1$ be a bounded domain with
smooth boundary. We are using in this paper the standard notation of function
spaces. By $L^{p}(\Omega)$ and $W^{k,p}(\Omega)$ with $1\leqslant
p\leqslant\infty$, $k\geqslant1$ we denote the Lebesgue spaces and
respectively, Sobolev spaces of functions on $\Omega$. If $X$ is a Banach
space with the norm $\Vert\cdot\Vert_{X}$, for $T>0$ we denote by
$L^{p}(0,T;X)$ the Banach space of all Bochner measurable functions
$u:(0,T)\rightarrow X$ such that $\Vert u\Vert_{X}\in L^{p}(0,T)$. \

Given a positive number $\nu$, we denote by $C^{\nu}(\overline{\Omega})$ the
H\"{o}lder space of $\left[  \nu\right]  $ times continuously differentiable
functions on $\overline{\Omega}$. We denote by $C^{\nu,\nu/2}(\overline
{\Omega}\times\left(  0,T\right)  )$ the H\"{o}lder space of exponents $\nu$
and $\nu/2$ by respect to $x$, respectively $t$ of continuous and bounded
functions defined on $\overline{\Omega}\times\left(  0,T\right)  $. If $J$ is
an interval of real numbers, the notation $C^{k}(J;X)$, $k\geqslant1$,
$k\in\mathbb{N}$ stands for the space of $k$ times continuously differentiable
functions from $J$ to the Banach space $X$.

Throughout this paper we denote by $C$, $C_{i}$ ($i=1,2,...$) positive
constants which may vary from line to line. These positive constants will be
independent of time, but we shall indicate explicitly on which other
parameters they are dependent, if it will be the case.\

Let $p\in(1,\infty)$ and $a$, $b$ be two positive constants. We denote
\[
A:=-a\Delta+b
\]
the positive self-adjoint operator with the domain defined by
\begin{equation}
\mathcal{D}(A):=\left\{  u\in W^{2,p}(\Omega);\;\;\frac{\partial u}{\partial
n}=0\quad\text{on\quad}\partial\Omega\right\}  . \label{dom}%
\end{equation}

For $0\leqslant\theta\leqslant1$ we denote the fractional powers of the
operator $A$ by $A^{\theta}:X_{p}^{\theta}\rightarrow L^{p}(\Omega)$ where the
space $X_{p}^{\theta}$ is endowed with the graph norm%
\[
\left\Vert u\right\Vert _{X_{p}^{\theta}}=\left\Vert A^{\theta}u\right\Vert
_{L^{p}(\Omega)}.
\]
Let us mention that $X_{p}^{1}=\mathcal{D}(A)$, $X_{p}^{0}=L^{p}(\Omega)$ and
$X_{p}^{\theta_{1}}$ is continuously embedded into $X_{p}^{\theta_{2}}$ if
$\theta_{1}\geqslant\theta_{2}$, moreover this embedding is compact if
$\theta_{1}>\theta_{2}$. We recall some results that we shall use throughout
the paper:

$(i)$ we have the embedding properties (see \cite[Theorem 1.6.1]{henry})
\begin{align}
X_{p}^{\theta}  &  \hookrightarrow W^{k,q}(\Omega), & k  &  -\frac{N}%
{q}<2\theta-\frac{N}{p},\text{\quad}q\geqslant p\label{h1}\\
X_{p}^{\theta}  &  \hookrightarrow C^{\nu}(\overline{\Omega}), & 0  &
\leqslant\nu<2\theta-\frac{N}{p}\,; \label{h2}%
\end{align}

$(ii)$ for all $u\in L^{p}(\Omega),$ $p\in(1,\infty)$ and $t>0$ there exists a
positive constant $C\left(  \theta\right)  $ such that (see \cite[Theorem
1.4.3]{henry})%
\begin{equation}
\Vert A^{\theta}e^{-tA}u\Vert_{L^{p}(\Omega)}\leqslant C\left(  \theta\right)
t^{-\theta}e^{-\delta t}\Vert u\Vert_{L^{p}(\Omega)},\quad\theta\geqslant0
\label{h3}%
\end{equation}
for some$\;\delta\in\left(  0,1\right)  $;

$(iii)$ for all $u\in L^{p}(\Omega),$ $1\leqslant q<p<\infty$ and $t>0$ we
have
\begin{equation}
\Vert A^{\beta}e^{-tA}u\Vert_{L^{p}(\Omega)}\leqslant C\left(  \beta\right)
t^{-\beta-\frac{N}{2}\left(  \frac{1}{q}-\frac{1}{p}\right)  }e^{-\delta
t}\Vert u\Vert_{L^{q}(\Omega)},\quad\beta\geqslant0 \label{h4}%
\end{equation}
for $\;\delta\in\left(  0,1\right)  $.

Throughout this paper we assume that $\chi$ and $g$ are non-negative functions
which satisfy the following conditions%
\begin{align*}
(H_{1})\qquad\chi &  \in C^{1}(\mathbb{R}_{+}),\;\chi\geqslant0\text{, }%
\chi\text{ and }\chi^{\prime}\text{ are globally Lipschitz continuous with
Lipschitz}\\
&  \text{constants }L_{\chi}\text{ and }L_{\chi^{\prime}}\text{ respectively;}%
\\
(H_{2})\qquad g  &  \in C^{1}(\mathbb{R}_{+}),\text{ }g\geqslant0\text{,
}g\text{ and }g^{\prime}\text{ are globally Lipschitz continuous with
Lipschitz}\\
&  \text{constants \ }L_{g}\text{ and }L_{g^{\prime}}\text{ respectively.}%
\end{align*}

\section{Local existence and uniqueness}

\label{local}

\setcounter{equation}{0} \setcounter{figure}{0}

In this section we establish the existence of local in time non-negative
solutions using a standard fixed point argument.

Let $p,q\in(1,\infty)$. We define $A_{1}=-\Delta+I$ and $A_{2}=-d\Delta+\gamma
I$ the positive definite self-adjoint operators with the domains
$\mathcal{D}(A_{1})=X_{q}^{1}$ and respectively $\mathcal{D}(A_{2})=X_{p}^{1}$
given by (\ref{dom}). Given $\tau>0$, $1<p<\infty$ and $\theta\in\left(
0,1\right)  $ we denote%
\[
Y_{q}=C\left(  [0,\tau];W^{1,q}(\Omega)\right)  ,\quad Z=C\left(
[0,\tau];X_{p}^{\theta}\right)  .
\]
We consider the closed set
\[
B_{\rho}^{\tau,q}:=\{(w,v,m)\in Y_{q}\times Y_{\infty}\times Z;\;\left\Vert
(w,v,m)\right\Vert _{Y_{q}\times Y_{\infty}\times Z}\leqslant\rho\}
\]
where $\left\Vert (w,v,m)\right\Vert _{Y_{q}\times Y_{\infty}\times Z}:=\Vert
w\Vert_{Y_{q}}+\Vert v\Vert_{Y_{\infty}}+\Vert m\Vert_{Z}$ and $\rho$\ is a
positive constant to be fixed later.

For $t\in\lbrack0,\tau]$ and $(w,v,m)\in B_{\rho}^{\tau,q}$ fixed we define
the mapping $F:=\left(  F_{1},F_{2},F_{3}\right)  $ by%
\begin{align*}
&  F_{1}(w,v,m)(t):=e^{-tA_{1}}w_{0}+\int\limits_{0}^{t}e^{-\left(
t-s\right)  A_{1}}G_{1}(w,v,m)(s)ds,\\
&  F_{2}(v,m)(t):=v_{0}-\int\limits_{0}^{t}G_{2}(v,m)(s)ds,\\
&  F_{3}(w,v)(t):=e^{-tA_{2}}m_{0}+\int\limits_{0}^{t}e^{-\left(  t-s\right)
A_{2}}G_{3}(w,v)(s)ds,
\end{align*}
where we denoted
\begin{align*}
&  G_{1}(w,v,m):=\chi(v)\nabla v\cdot\nabla w+\left(  \mu+1\right)  w-\mu
w\left(  wz^{-1}+v\right)  +\chi(v)wvm,\\
&  G_{2}(v,m):=mv,\\
&  G_{3}(w,v):=wz^{-1}g(v).
\end{align*}

\begin{lemma}
\label{fixed_point}Let $\Omega\subset\mathbb{R}^{N}$, $N\geqslant1$ be a
domain with $C^{2}$ boundary and $p>N$. Given an initial value $(w_{0}%
,v_{0},m_{0})\in W^{1,q}(\Omega)\times W^{1,\infty}(\Omega)\times
X_{p}^{\theta}$, where $\theta\in\left(  \frac{N+p}{2p},1\right]  $ and
$q\geqslant\frac{Np}{N+p}$ if $N\geqslant2$, there exists $\tau_{0}>0$
(depending only on $\left\Vert u_{0}\right\Vert _{W^{1,q}(\Omega)}$,
$\left\Vert v_{0}\right\Vert _{W^{1,\infty}(\Omega)}$ and $\left\Vert
m_{0}\right\Vert _{X_{p}^{\theta}}$) such that, for all $\tau\in(0,\tau_{0}]$
the application $F$ is a contraction from $B_{\rho}^{\tau,q}$\ into itself.
\end{lemma}

\begin{proof}
We shall prove first that the closed set $B_{\rho}^{\tau,q}$ is invariant by
$F$ for all $\tau\in(0,\tau_{0}]$ where $\tau_{0}$\ will be chosen later.
Taking $(w,v,m)\in B_{\rho}^{\tau,q}$, using the embeddings (\ref{h1}),
(\ref{h2}), the estimate (\ref{h3}) and the hypothesis $(H1)$ we estimate%
\begin{align}
&  \left\Vert F_{1}(w,v,m)(t)\right\Vert _{W^{1,q}(\Omega)}\leqslant\left\Vert
e^{-tA_{1}}w_{0}\right\Vert _{W^{1,q}(\Omega)}+C(\theta)\int\limits_{0}%
^{t}\left\Vert A_{1}^{\theta}e^{-\left(  t-s\right)  A_{1}}G_{1}%
(w,v,m)\right\Vert _{L^{q}\left(  \Omega\right)  }ds\leqslant\nonumber\\
&  \leqslant C_{1}\left\Vert w_{0}\right\Vert _{W^{1,q}(\Omega)}+C(\theta
)\int\limits_{0}^{t}\left(  t-s\right)  ^{-\theta}\left\Vert G_{1}%
(w,v,m)\right\Vert _{L^{q}\left(  \Omega\right)  }ds\leqslant\nonumber\\
&  \leqslant C_{1}\left\Vert w_{0}\right\Vert _{W^{1,q}(\Omega)}+C(\theta
)\rho\left[  \rho\left(  \rho+1\right)  \left(  L_{\chi}\rho+\chi(0)\right)
+\mu\rho\left(  1+C_{2}\left(  \rho\right)  \right)  +\left(  \mu+1\right)
\right]  \tau^{1-\theta},\quad\forall t\in\left[  0,\tau\right]  \label{f1}%
\end{align}
where $C_{2}\left(  \rho\right)  =e^{\rho\left[  2^{-1}L_{\chi}\rho
+\chi(0)\right]  }$. We also have%
\begin{equation}
\Vert F_{2}(v,m)(t)\Vert_{W^{1,\infty}(\Omega)}\leqslant\Vert v_{0}%
\Vert_{W^{1,\infty}(\Omega)}+2\rho^{2}\tau,\qquad\forall t\in\left[
0,\tau\right]  . \label{f2}%
\end{equation}
Taking into account the embeddings (\ref{h1}), (\ref{h2}), the estimate
(\ref{h3}) and the hypothesis $(H2)$, we obtain
\begin{align}
&  \Vert F_{3}(w,v)(t)\Vert_{X_{p}^{\theta}}\leqslant\left\Vert e^{-tA_{2}%
}m_{0}\right\Vert _{X_{p}^{\theta}}+C\left(  \theta\right)  \int
\limits_{0}^{t}\left\Vert A_{2}^{\theta}e^{-\left(  t-s\right)  A_{2}}%
wz^{-1}g(v)\right\Vert _{L^{p}\left(  \Omega\right)  }ds\leqslant\nonumber\\
&  \leqslant C_{3}\left\Vert m_{0}\right\Vert _{X_{p}^{\theta}}+C\left(
\theta\right)  \int\limits_{0}^{t}\left(  t-s\right)  ^{-\theta}\left\Vert
wz^{-1}g(v)\right\Vert _{L^{p}\left(  \Omega\right)  }ds\leqslant\nonumber\\
&  \leqslant C_{3}\left\Vert m_{0}\right\Vert _{X_{p}^{\theta}}+C\left(
\theta\right)  C_{2}\left(  \rho\right)  \rho\left(  L_{g}\rho+g(0)\right)
\tau^{1-\theta},\qquad\forall t\in\left[  0,\tau\right]  . \label{f3}%
\end{align}
Denoting $C_{4}=\max\left\{  1,C_{1},C_{3}\right\}  $, the estimates
(\ref{f1}), (\ref{f2}) and (\ref{f3}) imply that there exists the constant
\[
C_{5}\left(  \rho\right)  =C(\theta)\rho\left[  \rho\left(  \rho+1\right)
\left(  L_{\chi}\rho+\chi(0)\right)  +C_{2}\left(  \rho\right)  \left(
L_{g}\rho+g(0)\right)  +\mu\rho\left(  1+C_{2}\left(  \rho\right)  \right)
+\left(  \mu+1\right)  \right]
\]
such that
\begin{align*}
&  \Vert F_{1}(w,v,m)\Vert_{Y_{q}}+\Vert F_{2}(v,m)\Vert_{Y_{\infty}}+\Vert
F_{3}(w,v)\Vert_{Z}\leqslant\\
&  \leqslant C_{4}\left(  \left\Vert w_{0}\right\Vert _{W^{1,q}(\Omega)}+\Vert
v_{0}\Vert_{W^{1,\infty}(\Omega)}+\left\Vert m_{0}\right\Vert _{X_{p}^{\theta
}}\right)  +C_{5}\left(  \rho\right)  \tau^{1-\theta}+2\rho^{2}\tau\,
\end{align*}
provided that $(w,v,m)\in B_{\rho}^{\tau,q}$. Now we fix $\rho>2C_{4}\left(
\left\Vert w_{0}\right\Vert _{W^{1,q}(\Omega)}+\Vert v_{0}\Vert_{W^{1,\infty
}(\Omega)}+\left\Vert m_{0}\right\Vert _{X_{p}^{\theta}}\right)  >0$
sufficiently large and we choose $\tau_{1}>0$ small enough such that
$F(B_{\rho}^{\tau,q})\subset B_{\rho}^{\tau,q}$ for all $\tau\in(0,\tau_{1}]$.

In what follows we prove that $F$ is a contraction. Let $(w,v,m),$
$(\overline{w},\overline{v},\overline{m})\in B_{\rho}^{\tau,q}$, where
$\tau\in(0,\tau_{1}]$ and $\rho$ was fixed previously. Taking into account
(\ref{h1}), (\ref{h2}), (\ref{h3}) and the hypothesis $(H1)$ we estimate%
\begin{align}
&  \Vert F_{1}(w,v,m)(t)-F_{1}(\overline{w},\overline{v},\overline{m}%
)(t)\Vert_{W^{1,q}(\Omega)}\leqslant\nonumber\\
&  \leqslant C(\theta)\int\limits_{0}^{t}\left\Vert A_{1}^{\theta}e^{-\left(
t-s\right)  A_{1}}(\chi(v)\nabla v\cdot\nabla w-\chi(\overline{v}%
)\nabla\overline{w}\cdot\nabla\overline{v})\right\Vert _{L^{q}\left(
\Omega\right)  }ds+\nonumber\\
&  +\left(  \mu+1\right)  C(\theta)\int\limits_{0}^{t}\left\Vert A_{1}%
^{\theta}e^{-\left(  t-s\right)  A_{1}}\left(  w-\overline{w}\right)
\right\Vert _{L^{q}\left(  \Omega\right)  }ds+\mu C(\theta)\int\limits_{0}%
^{t}\left\Vert A_{1}^{\theta}e^{-\left(  t-s\right)  A_{1}}\left[  \left(
w^{2}z^{-1}-\overline{w}^{2}\overline{z}^{-1}\right)  \right]  \right\Vert
_{L^{q}\left(  \Omega\right)  }ds+\nonumber\\
&  +\mu C(\theta)\int\limits_{0}^{t}\left\Vert A_{1}^{\theta}e^{-\left(
t-s\right)  A_{1}}\left(  wv-\overline{w}\overline{v}\right)  \right\Vert
_{L^{q}\left(  \Omega\right)  }ds+\mu C(\theta)\int\limits_{0}^{t}\left\Vert
A_{1}^{\theta}e^{-\left(  t-s\right)  A_{1}}\left(  \chi(v)wvm-\chi
(\overline{v})\overline{w}\overline{v}\overline{m}\right)  \right\Vert
_{L^{q}\left(  \Omega\right)  }ds\leqslant\nonumber\\
&  \leqslant C_{6}(\rho)\tau^{1-\theta}\left[  \left\Vert w-\overline
{w}\right\Vert _{Y_{q}}+\left\Vert v-\overline{v}\right\Vert _{Y_{\infty}%
}+\left\Vert m-\overline{m}\right\Vert _{Z}\right]  ,\qquad\forall t\in\left[
0,\tau\right]  \label{ff1}%
\end{align}
where%
\[
C_{6}(\rho)=C(\theta)\left\{  \rho\left(  L_{\chi}\rho+\chi(0)\right)  \left[
1+\mu\rho\left(  1+2C_{2}(\rho)\right)  \right]  +2\mu\rho C_{2}(\rho)+\left(
\mu+1\right)  \right\}  .
\]
Also we obtain%
\begin{equation}
\Vert F_{2}(w,v,m)(t)-F_{2}(\overline{w},\overline{v},\overline{m}%
)(t)\Vert_{W^{1,\infty}(\Omega)}\leqslant2\rho\tau\left(  \left\Vert
v-\overline{v}\right\Vert _{Y_{\infty}}+\left\Vert m-\overline{m}\right\Vert
_{Z}\right)  ,\qquad\forall t\in\left[  0,\tau\right]  . \label{ff2}%
\end{equation}
Using (\ref{h1}), (\ref{h2}), (\ref{h3}) and the hypothesis $(H2)$ we get
\begin{align}
&  \Vert F_{3}(w,v)(t)-F_{3}(\overline{w},\overline{v})(t)\Vert_{X_{p}%
^{\theta}}\leqslant C\left(  \theta\right)  \int\limits_{0}^{t}\left\Vert
A_{2}^{\theta}e^{-\left(  t-s\right)  A_{2}}\left[  wz^{-1}g(v)-\overline
{w}\overline{z}^{-1}g(\overline{v})\right]  \right\Vert _{L^{p}\left(
\Omega\right)  }ds\leqslant\nonumber\\
&  \leqslant C_{7}(\rho)\tau^{1-\theta}\left[  \left\Vert w-\overline
{w}\right\Vert _{Y_{q}}+\left\Vert v-\overline{v}\right\Vert _{Y_{\infty}%
}\right]  ,\qquad\forall t\in\left[  0,\tau\right]  \label{ff3}%
\end{align}
where%
\[
C_{7}(\rho)=C\left(  \theta\right)  C_{2}(\rho)\left(  L_{g}\rho+g(0)\right)
\left(  \rho+1\right)  .
\]

Thus, taking $\tau_{0}$ sufficiently small such that $\tau_{0}<\tau_{1}$, we
obtain from (\ref{ff1}),(\ref{ff2}) and (\ref{ff3}) that the mapping $F$ is a
contraction from $B_{\rho}^{\tau,q}$ into itself for all $\tau\in(0,\tau_{0}]$.
\end{proof}

We shall prove now the existence of a unique non-negative maximal solution to
(\ref{e1})-(\ref{i}).

\begin{theorem}
\label{local_existence}Let $\Omega\subset\mathbb{R}^{N}$, $N\geqslant1$ be a
domain with $C^{2}$ boundary and $p>N$. Given $(w_{0},v_{0},m_{0})\in
W^{1,q}(\Omega)\times W^{1,\infty}(\Omega)\times X_{p}^{\theta}$, $\theta
\in\left(  \frac{N+p}{2p},1\right)  $, $q\geqslant\frac{Np}{N+p}$ if
$N\geqslant2$ and $r=\max\left\{  2q,p\right\}  $, there exists $T>0$
(depending only on $\left\Vert w_{0}\right\Vert _{W^{1,q}(\Omega)}$,
$\left\Vert v_{0}\right\Vert _{W^{1,\infty}(\Omega)}$ and $\left\Vert
m_{0}\right\Vert _{X_{p}^{\theta}}$) such that the problem (\ref{e1}%
)-(\ref{i}) has a unique solution $\left(  w,v,m\right)  $ defined on an
interval $[0,T)\subset\mathbb{R}$ and%
\begin{align*}
w  &  \in C\left(  \left[  0,T\right)  ;W^{1,q}(\Omega)\right)  \cap C\left(
\left(  0,T\right)  ;W^{2,r}(\Omega)\right)  \cap C^{1}\left(  \left(
0,T\right)  ;W^{1,q}(\Omega)\right) \\
v  &  \in C\left(  \left[  0,T\right)  ;W^{1,\infty}(\Omega)\right)  \cap
C^{1}\left(  (0,T);W^{1,\infty}(\Omega)\right) \\
m  &  \in C\left(  \left[  0,T\right)  ;X_{p}^{\theta}\right)  \cap C\left(
\left(  0,T\right)  ;W^{2,p}(\Omega)\right)  \cap C^{1}\left(  \left(
0,T\right)  ;X_{p}^{\theta}\right)  .
\end{align*}
Moreover, the solution depends continuously on the initial data.
\end{theorem}

\begin{proof}
Lemma \ref{fixed_point} shows that the map $F:B_{\rho}^{\tau,q}\rightarrow
B_{\rho}^{\tau,q}$, $\tau\in(0,\tau_{0}]$, has a unique fixed point $(w,v,m)$
which is the weak solution to the system (\ref{e1})-(\ref{i}). Taking into
account the fact that the map $F$ is a contraction and the estimates
established in Lemma \ref{fixed_point} we obtain the continuous dependence of
this solution on the initial data.

In what follows we show the existence of a unique maximal solution to
(\ref{e1})-(\ref{i}) having the regularity properties stated in the theorem.
First we can prove that for every fixed $t\in(0,\tau_{0}]$ the maps
$G_{1}(t):W^{1,q}(\Omega)\times W^{1,\infty}(\Omega)\times X_{p}^{\theta
}\rightarrow L^{q}(\Omega)$, $G_{2}(t):W^{1,\infty}(\Omega)\times
X_{p}^{\theta}\rightarrow W^{1,\infty}(\Omega)$ and $G_{3}(t):W^{1,q}%
(\Omega)\times W^{1,\infty}(\Omega)\rightarrow L^{p}(\Omega)$ are Lipschitzian
using similar arguments to those used in the proof of Lemma \ref{fixed_point}.
Therefore applying \cite[Theorem 3.5.2]{henry} we obtain%
\begin{equation}
w\in C^{1}\left(  (0,\tau_{0}];W^{1,q}(\Omega)\right)  ,\quad m\in
C^{1}\left(  (0,\tau_{0}];X_{p}^{\theta}\right)  . \label{wm}%
\end{equation}
Let us observe that for every $t\in(0,\tau_{0}],$ $w(t)\in W^{1,q}(\Omega)$ is
the solution to the problem%
\begin{align*}
&  -\Delta w(t)-a(t)\cdot\nabla w(t)=f(t)-\frac{\partial w}{\partial t}(t) &
x  &  \in\Omega\\
&  \frac{\partial w}{\partial\eta}=0 & x  &  \in\partial\Omega
\end{align*}
where we have denoted%
\[
a(t):=\chi(v(t))\nabla v(t)\in\left(  L^{\infty}(\Omega)\right)  ^{N},\quad
f(t):=\left(  \mu w\left(  1-wz^{-1}-v\right)  +\chi(v)wvm\right)  (t)\in
L^{q}(\Omega).
\]
The elliptic regularity implies that $w(t)\in W^{2,q}(\Omega)$ for every
$t\in(0,\tau_{0}]$. As $q>N/2$, we deduce from the Sobolev embeddings that
$w(t)\in W^{2,2q}(\Omega)$. Finally, by recurrence, we obtain $w(t)\in
W^{2,r}(\Omega)$, $r=\max\left\{  2q,p\right\}  $. In a similar manner we can
prove%
\[
m\in C\left(  (0,\tau_{0}];W^{2,p}(\Omega)\right)  ,\quad2\leqslant p<\infty.
\]

In order to show the uniqueness of the solutions in the spaces indicated
above, let us suppose that $\left(  w_{1},v_{1},m_{1}\right)  $, $\left(
w_{2},v_{2},m_{2}\right)  $ are two different solutions to the system
(\ref{e1})-(\ref{i}). Let $\left[  0,T\right)  $ the maximal interval where
both solutions are defined and
\[
I=\left\{  t\in\left[  0,T\right)  ;\quad w_{1}(x,t)=w_{2}(x,t),\ v_{1}%
(x,t)=v_{2}(x,t),\ m_{1}(x,t)=m_{2}(x,t),\ \forall x\in\Omega\right\}  .
\]
From the local existence we deduce that $I$ is a nonempty set, so let
$\tau_{\max}>0$ be maximal such that $\left[  0,\tau_{\max}\right)  \subset I$
and suppose $\tau_{\max}<T$. As $I$ is a closed set in $\left[  0,T\right)  $
(from the continuity of the solutions), it follows that $\left[  0,\tau_{\max
}\right]  \subset I$. Applying now the local existence result with the initial
conditions considered in $\tau_{\max}$ it follows that the system
(\ref{e1})-(\ref{i}) has a unique solution on a small interval $\left(
\tau_{\max},\tau_{\max}+\varepsilon\right)  $. This contradicts the maximality
of $\tau_{\max}$, hence $\tau_{\max}=T$.

Let us remark that the choice of $\tau_{0}$ in Lemma \ref{fixed_point} depends
only on the initial data, respectively $\left\Vert w_{0}\right\Vert
_{W^{1,q}(\Omega)}$, $\left\Vert v_{0}\right\Vert _{W^{1,p}(\Omega)}$ and
$\left\Vert m_{0}\right\Vert _{X_{p}^{\theta}}$, so the solution $(u,v,m)$ can
be extended up to a maximal time $\tau_{\max}=T$ that depends also only on the
initial data.
\end{proof}

\begin{remark}
Let us mention that in fact we can obtain from the Theorem
\ref{local_existence} a better regularity for $w$. Indeed, $G_{1}(t)\in
C^{\alpha}(\Omega)$ for some $0<\alpha<1$ and $w(t,\cdot)\in C^{2+\alpha
}(\Omega)$. Thus for $t>0$, $\left(  t,x\right)  \rightarrow w(t,x;w_{0})$ is
continuously differentiable in $t$, twice continuously differentiable in $x$,
and hence $w$ is a classical solution.
\end{remark}

In view of (\ref{changevariables2}) we observe that Theorem
\ref{local_existence} implies also the local existence of the solution to the
initial system (\ref{eq1})-(\ref{ic}). Moreover, we shall show the
non-negativity of local solutions to (\ref{eq1})-(\ref{ic}) corresponding to
non-negative initial values $\left(  w_{0},v_{0},m_{0}\right)  $.

\begin{theorem}
\label{local_existence_u}Let $\Omega\subset\mathbb{R}^{N}$, $N\geqslant1$ be a
domain with $C^{2}$ boundary and $p>N$. Given the non-negative initial value
$(u_{0},v_{0},m_{0})\in W^{1,p}(\Omega)\times W^{1,\infty}(\Omega)\times
X_{p}^{\theta}$, $\theta\in\left(  \frac{N+p}{2p},1\right)  $, there exists
$T>0$ (depending only on $\left\Vert w_{0}\right\Vert _{W^{1,p}(\Omega)}$,
$\left\Vert v_{0}\right\Vert _{W^{1,\infty}(\Omega)}$ and $\left\Vert
m_{0}\right\Vert _{X_{p}^{\theta}}$) such that the problem (\ref{eq1}%
)-(\ref{ic}) has a unique non-negative solution $\left(  u,v,m\right)  $
defined on an interval $[0,T)\subset\mathbb{R}$ and%
\begin{align*}
u  &  \in C\left(  \left[  0,T\right)  ;W^{1,q}(\Omega)\right)  \cap C\left(
\left(  0,T\right)  ;W^{1,\infty}(\Omega)\right)  \cap C^{1}\left(  \left(
0,T\right)  ;W^{1,q}(\Omega)\right) \\
v  &  \in C\left(  \left[  0,T\right)  ;W^{1,\infty}(\Omega)\right)  \cap
C^{1}\left(  (0,T);W^{1,\infty}(\Omega)\right) \\
m  &  \in C\left(  \left[  0,T\right)  ;X_{p}^{\theta}\right)  \cap C\left(
\left(  0,T\right)  ;W^{2,p}(\Omega)\right)  \cap C^{1}\left(  \left(
0,T\right)  ;X_{p}^{\theta}\right)  .
\end{align*}
Moreover, the solution depends continuously on the initial data.
\end{theorem}

\begin{proof}
From (\ref{changevariables2}) and taking into account $(u_{0},v_{0})\in
W^{1,p}(\Omega)\times W^{1,\infty}(\Omega)$ we obtain $w_{0}\in W^{1,p}%
(\Omega)$ and the hypotheses of Theorem \ref{local_existence}\ are satisfied.
Using again (\ref{changevariables2}) we recover the regularity of $u$.

Let us observe that from the equation (\ref{eq2}) we obtain
\begin{equation}
v=v_{0}e^{-\int\limits_{0}^{t}m(x,s)ds} \label{v}%
\end{equation}
which implies the nonnegativity of $v$. We shall prove the nonnegativity of
the solution $u$ by the truncation technique used in \cite[Theorem 2.1]{yagi}.
According to \cite{yagi}, there exists a decreasing function $H\in
\mathcal{C}^{3}(\mathbb{R};\overline{\mathbb{R}}_{+})$ and a constant
$C_{56}>0$ such that $H(u)>0$ if $u<0$ and $H(u)=0$ if $u\geqslant0$ and
having the properties%
\begin{align}
0  &  \leqslant H^{\prime\prime}(u)u^{2}\leqslant\widetilde{C}H(u),\qquad
u\in\mathbb{R}\label{H1}\\
0  &  \leqslant H^{\prime}(u)u\leqslant\widetilde{C}H(u),\qquad\ \ u\in
\mathbb{R}\label{H2}\\
0  &  \leqslant H(u)\leqslant\widetilde{C}u^{2},\qquad\ \ \ \ \ \ \ \ \ u\in
\mathbb{R}. \label{H4}%
\end{align}
We consider the non-negative function%
\[
\varphi(t)=%
{\displaystyle\int\limits_{\Omega}}
H(u(x,t))dx,\qquad0\leqslant t\leqslant\tau.
\]
The definition of the function $H(u)$ implies that $\varphi\in\mathcal{C}%
(\left[  0,\tau\right]  ;\overline{\mathbb{R}}_{+})\cap\mathcal{C}^{1}%
((0,\tau];\mathbb{R})$, $\varphi(0)=0$ and $\varphi$ has the derivative%
\begin{align*}
\varphi^{\prime}(t)  &  =%
{\displaystyle\int\limits_{\Omega}}
H^{\prime}(u)\frac{\partial u}{\partial t}dx=-%
{\displaystyle\int\limits_{\Omega}}
H^{\prime\prime}(u)\left\vert \nabla u\right\vert ^{2}dx+%
{\displaystyle\int\limits_{\Omega}}
u\chi(v)H^{\prime\prime}(u)\nabla u\cdot\nabla vdx+\mu%
{\displaystyle\int\limits_{\Omega}}
H^{\prime}(u)u(1-u-v)dx\leqslant\\
&  \leqslant\frac{1}{2}\widetilde{C}%
{\displaystyle\int\limits_{\Omega}}
H(u)\chi^{2}(v)\left\vert \nabla v\right\vert ^{2}dx+\mu\widetilde{C}%
{\displaystyle\int\limits_{\Omega}}
H(u)dx\leqslant C_{8}\varphi(t)\qquad
\end{align*}
for all $0\leqslant t\leqslant\tau$\ where%
\[
C_{8}=\widetilde{C}\left[  \frac{1}{2}\left(  L_{\chi}\left\Vert v\right\Vert
_{L^{\infty}(\Omega)}+\chi(0)\right)  ^{2}\left\Vert \nabla v\right\Vert
_{L^{\infty}(\Omega)}^{2}+\mu\right]  .
\]
Thus the Gronwall inequality ensures $\varphi(t)=0$ for all $t\in\left[
0,\tau\right]  $, that is $u(t)\geqslant0$ on $\Omega$ for $t\in\left[
0,\tau\right]  $. Finally, since $ug(v)\ $is a non-negative function it is
straightforward to prove that $m(t)\geqslant0$ on $\Omega$ for $t\in\left[
0,\tau\right]  $ using the maximum principle for parabolic equations. But
$\tau>0$ is arbitrary, so the desired positivity follows.
\end{proof}

\begin{theorem}
\label{local_existence copy_u2}Under the same hypotheses as in Theorem
\ref{local_existence_u}, if $v_{0}\in W^{2,r}(\Omega)$, $r=\max\left\{
2q,p\right\}  $ then the problem (\ref{eq1})-(\ref{ic}) has a unique
non-negative solution $\left(  u,v,m\right)  $ defined on an interval
$[0,T)\subset\mathbb{R}$ and%
\begin{align*}
u  &  \in C\left(  \left[  0,T\right)  ;W^{1,q}(\Omega)\right)  \cap C\left(
\left(  0,T\right)  ;W^{2,r}(\Omega)\right)  \cap C^{1}\left(  \left(
0,T\right)  ;W^{1,q}(\Omega)\right) \\
v  &  \in C\left(  \left[  0,T\right)  ;W^{2,r}(\Omega)\right)  \cap
C^{1}\left(  (0,T);W^{2,r}(\Omega)\right) \\
m  &  \in C\left(  \left[  0,T\right)  ;X_{p}^{\theta}\right)  \cap C\left(
\left(  0,T\right)  ;W^{2,p}(\Omega)\right)  \cap C^{1}\left(  \left(
0,T\right)  ;X_{p}^{\theta}\right)  .
\end{align*}
Moreover, the solution depends continuously on the initial data.
\end{theorem}

\section{Global existence in time}

\label{global}\

\setcounter{equation}{0} \setcounter{figure}{0}

We denote by $(u,v,m)$ the maximal non-negative local solution to the problem
(\ref{eq1})-(\ref{ic}) on $[0,T)$ and we shall prove that $T=+\infty$ in the
case $N=3$. Throughout this section all the constants are independent of $T$
and when it will be the case we shall make explicit their dependence on the
data of the problem.

\begin{lemma}
\label{mass} Let $\Omega\subset\mathbb{R}^{N}$, $N\geqslant1$ be a domain with
smooth boundary and $\left(  u_{0},v_{0},m_{0}\right)  \in L^{1}(\Omega
)\times$ $L^{\infty}(\Omega)\times$ $L^{1}(\Omega)$. Then the solution
$(u,v,m)$ to the problem (\ref{eq1})-(\ref{ic}) satisfies the following
estimates
\begin{align}
&  \Vert u(\cdot,t)\Vert_{L^{1}(\Omega)}\leqslant\max\{|\Omega|,\Vert
u_{0}\Vert_{L^{1}(\Omega)}\}\;,\;\;\;\label{mass1}\\
&  \Vert v(\cdot,t)\Vert_{L^{\infty}(\Omega)}\leqslant\Vert v_{0}%
\Vert_{L^{\infty}(\Omega)},\;\;\;\label{mass2}\\
&  \Vert m(\cdot,t)\Vert_{L^{1}(\Omega)}\leqslant\Vert m_{0}\Vert
_{L^{1}(\Omega)}e^{-t}+\left[  L_{g}\Vert v_{0}\Vert_{L^{\infty}(\Omega
)}+g(0)\right]  \max\{|\Omega|,\Vert u_{0}\Vert_{L^{1}(\Omega)}\}
\label{mass3}%
\end{align}
for all $t>0$.
\end{lemma}

\begin{proof}
Integrating the equation (\ref{eq1}) in space and taking into account the
non-negativity of the solution and the boundary condition (\ref{bc}) we get
\[
\frac{d}{dt}\int_{\Omega}udx\leqslant\mu\int_{\Omega}u(x,t)dx-\mu\int_{\Omega
}u^{2}(x,t)dx.
\]
Applying Jensen's inequality and Gronwall's lemma we obtain (\ref{mass1}).

The boundedness of $v$ results immediately from the equation (\ref{eq2}). We
obtain (\ref{mass3}) integrating the equation (\ref{eq3}) in space and taking
into account the boundary condition (\ref{bc}), the estimates (\ref{mass1}),
(\ref{mass2}) and the Gronwall lemma.
\end{proof}

\begin{corollary}
\label{mass w} Let $\Omega\subset\mathbb{R}^{N}$, $N\geqslant1$ be a domain
with smooth boundary and $u_{0}\in L^{1}(\Omega)$ (or, equivalently $\left(
w_{0},v_{0}\right)  \in L^{1}(\Omega)\times L^{\infty}(\Omega)$). Then
\begin{equation}
\Vert w(\cdot,t)\Vert_{L^{1}(\Omega)}\leqslant\max\{|\Omega|,\Vert u_{0}%
\Vert_{L^{1}(\Omega)}\}\;,\;\;\; \label{mass4}%
\end{equation}
for all $t>0$, where $w(x,t)$ is the function given by (\ref{changevariables2}).
\end{corollary}

In what follows we show that if the initial data $m_{0}$ is in an appropriate
space, then we can find a bound for $m$ in $L^{j}(\Omega)$, $j>1$ for $t>0$.
This result is based on \cite[Lemma 4.1]{horstmann}.

\begin{lemma}
\label{h41}Let $\Omega\subset\mathbb{R}^{N}$, $N\geqslant2$ be a domain with
smooth boundary. Assume that there exist $r\in\left[  1,N\right)  $ and
$1<p<rN/(N-r)$ such that%
\[
\left\Vert u(\cdot,t)\right\Vert _{L^{r}(\Omega)}\leqslant C_{9}%
\]
for all $t\in\left(  0,T\right)  $ and $\left(  v_{0},m_{0}\right)  \in
L^{\infty}(\Omega)\times$ $W^{1,p}(\Omega)$. Then%
\begin{equation}
\left\Vert m\left(  \cdot,t\right)  \right\Vert _{W^{1,p}}\leqslant C\left(
p,r,\Vert v_{0}\Vert_{L^{\infty}(\Omega)},\Vert m_{0}\Vert_{W^{1,p}(\Omega
)}\right)  \left(  1+C_{9}\right)  \label{mass5}%
\end{equation}
for all $t\in\left(  0,T\right)  $.
\end{lemma}

\begin{proof}
We fix $1<p<rN/(N-r)$ and we choose $\beta$ such that
\[
\frac{1}{2}<\beta<\frac{1}{2}+\frac{N}{2}\left(  \frac{1}{p}-\frac{N-r}%
{rN}\right)  .
\]
From the representation formula
\begin{equation}
m(t)=e^{-tA_{2}}m_{0}+\int\limits_{0}^{t}e^{-\left(  t-s\right)  A_{2}%
}u(s)g(v(s))ds,\quad\forall t\in(0,T) \label{rf}%
\end{equation}
and taking into account (\ref{h4}) we obtain
\begin{align*}
&  \left\Vert m\left(  \cdot,t\right)  \right\Vert _{W^{1,p}}\leqslant
C\left(  \beta\right)  \Vert m_{0}\Vert_{W^{1,p}}+\int\limits_{0}%
^{t}\left\Vert A_{2}^{\beta}e^{-\left(  t-s\right)  A_{2}}%
u(s)g(v(s)\right\Vert _{L^{p}(\Omega)}ds\leqslant\\
&  \leqslant C\left(  \beta\right)  \Vert m_{0}\Vert_{W^{1,p}}+\frac{C\left(
\beta,p\right)  C_{9}}{^{1-\beta-\frac{N}{2}\left(  \frac{1}{r}-\frac{1}%
{p}\right)  }}\left[  L_{g}\Vert v_{0}\Vert_{L^{\infty}(\Omega)}+g(0)\right]
,\quad\forall t\in(0,T)
\end{align*}
Using Lemma \ref{mass} the statement follows.
\end{proof}

\begin{remark}
If $m_{0}\in L^{1}(\Omega)$, the estimate (\ref{mass5}) is still valid for all
$t\in\left[  \tau,T\right)  $ where $\tau\in\left(  0,\min\left\{
1,T\right\}  \right)  $.
\end{remark}

\begin{corollary}
\label{cor_h1}Let $\Omega\subset\mathbb{R}^{N}$, $N>2$ be a domain with smooth
boundary and assume that the hypotheses of Lemma \ref{h41} are satisfied.
If\ $r\in\left[  1,N/2\right)  $ then%
\begin{equation}
\left\Vert m\left(  \cdot,t\right)  \right\Vert _{L^{j}}\leqslant C_{10}
\label{mass6}%
\end{equation}
for all $t\in(0,T)$ and $\frac{N}{N-1}<j<\frac{rN}{N-2r}$.
\end{corollary}

The next proposition asserts that if the initial value is in a suitable space
and (\ref{mass6}) is satisfied, then we obtain an estimate for $w$, and
respectively for $u$, in $L^{q}(\Omega)$, for all $q>1$ and $t>0$. Moreover,
we shall show later that if $w_{0}\in L^{\infty}(\Omega)$ then we derive a
bound for $w$ in $L^{\infty}(\Omega)$ for all $t>0.$

\begin{proposition}
\label{Lp}Let $\Omega\subset\mathbb{R}^{3}$ be a domain with smooth boundary
and $(w_{0},v_{0},m_{0})\in L^{q}(\Omega)\times L^{\infty}(\Omega)\times
W^{1,p}$, $1<q<\infty$, $p\geqslant\frac{6}{5}$. Then there exists a constant
independent on time $C_{11}=C_{11}(p,\left\Vert w_{0}\right\Vert
_{L^{q}(\Omega)},\left\Vert v_{0}\right\Vert _{L^{\infty}(\Omega)})$ such that
the solution $w(x,t)$ to the system (\ref{e1})-(\ref{i}) satisfies \
\begin{equation}
\left\Vert w\right\Vert _{L^{\infty}(0,T;L^{q}(\Omega))}\leqslant
C_{11},\qquad\forall~1<q<+\infty. \label{lp}%
\end{equation}

\end{proposition}

\begin{proof}
Multiplying the equation (\ref{e1}) with $qw^{q-1}z^{-1}$, $q\in
\lbrack2,\infty)$ and integrating in space we have
\begin{align}
\frac{d}{dt}%
{\displaystyle\int\limits_{\Omega}}
z^{-1}w^{q}  &  =-\frac{4(q-1)}{q}%
{\displaystyle\int\limits_{\Omega}}
z^{-1}\left\vert \nabla w^{q/2}\right\vert ^{2}+\mu q%
{\displaystyle\int\limits_{\Omega}}
z^{-1}w^{q}\left(  1-wz^{-1}-v\right)  +\nonumber\\
&  +\left(  q-1\right)
{\displaystyle\int\limits_{\Omega}}
z^{-1}\chi(v)w^{q}vm. \label{estw}%
\end{align}
On both sides of (\ref{estw}) we add the term $\varepsilon\Vert w^{q/2}%
\Vert_{L^{2}(\Omega)}^{2}$ where $\varepsilon>0$\ is a constant to be
determined later. Taking into account the non-negativity of the solution and
Lemma \ref{mass} we obtain%
\begin{align}
\frac{d}{dt}%
{\displaystyle\int\limits_{\Omega}}
z^{-1}w^{q}+\varepsilon\Vert w^{q/2}\Vert_{L^{2}(\Omega)}^{2}  &
\leqslant-\frac{4(q-1)}{q}%
{\displaystyle\int\limits_{\Omega}}
z^{-1}\left\vert \nabla w^{q/2}\right\vert ^{2}+\left(  q\mu C_{12}%
+\varepsilon\right)  \Vert w^{q/2}\Vert_{L^{2}(\Omega)}^{2}+\nonumber\\
&  +\left(  q-1\right)  C_{13}%
{\displaystyle\int\limits_{\Omega}}
w^{q}m, \label{es01}%
\end{align}
where%
\begin{align*}
C_{12}  &  =C_{12}\left(  \Vert v_{0}\Vert_{L^{\infty}(\Omega)}\right)
=e^{\frac{L_{\chi}}{2}\Vert v_{0}\Vert_{L^{\infty}(\Omega)}^{2}+\chi(0)\Vert
v_{0}\Vert_{L^{\infty}(\Omega)}},\\
C_{13}  &  =C_{13}\left(  \Vert v_{0}\Vert_{L^{\infty}(\Omega)}\right)
=e^{\frac{L_{\chi}}{2}\Vert v_{0}\Vert_{L^{\infty}(\Omega)}^{2}+\chi(0)\Vert
v_{0}\Vert_{L^{\infty}(\Omega)}}\left(  L_{\chi}\Vert v_{0}\Vert_{L^{\infty
}(\Omega)}+\chi(0)\right)  \Vert v_{0}\Vert_{L^{\infty}(\Omega)}.
\end{align*}
Now we estimate the last two terms from the right-hand side of (\ref{es01}).
Using Gagliardo-Nirenberg and Young's inequalities and Lemma \ref{mass} we
obtain
\begin{align}
\left(  q\mu C_{12}+\varepsilon\right)  \Vert w^{q/2}\Vert_{L^{2}(\Omega
)}^{2}  &  \leqslant C(\Omega)\left(  q\mu C_{12}+\varepsilon\right)  \Vert
w^{q/2}\Vert_{W^{1,2}(\Omega)}^{6\left(  q-1\right)  /\left(  3q-1\right)
}\Vert w\Vert_{L^{1}(\Omega)}^{2q/\left(  3q-1\right)  }\leqslant\nonumber\\
&  \leqslant\frac{\varepsilon}{4}\Vert w^{q/2}\Vert_{W^{1,2}(\Omega)}%
^{2}+C_{14}\Vert w\Vert_{L^{1}(\Omega)}^{4q/\left(  3q-1\right)  }, \label{t1}%
\end{align}
where%
\[
C_{14}=C_{14}(\varepsilon,q,\Vert v_{0}\Vert_{L^{\infty}(\Omega)},\left\vert
\Omega\right\vert )=\left(  \frac{2}{3\left(  q-1\right)  }\right)  \left(
\frac{\varepsilon}{4}\right)  ^{-3\left(  q-1\right)  /2}\left(
\frac{3\left(  q-1\right)  }{3q-1}C(\Omega)\left(  q\mu C_{12}+\varepsilon
\right)  \right)  ^{\left(  3q-1\right)  /2}.
\]
Taking into account H\"{o}lder, Gagliardo-Nirenberg and Young's inequalities,
the boundedness (\ref{mass6}) and Lemma \ref{mass}\ we estimate the last term
from the right-hand side of (\ref{es01})
\begin{align}
\left(  q-1\right)  C_{13}%
{\displaystyle\int\limits_{\Omega}}
w^{q}m  &  \leqslant\left(  q-1\right)  C_{13}\Vert w^{q/2}\Vert_{L^{4}%
(\Omega)}^{2}\Vert m\Vert_{L^{2}(\Omega)}\leqslant\nonumber\\
&  \leqslant\left(  q-1\right)  C_{13}C(\Omega)\Vert w^{q/2}\Vert
_{W^{1,2}(\Omega)}^{3\left(  2q-1\right)  /\left(  3q-1\right)  }\Vert
w^{q/2}\Vert_{L^{2/q}(\Omega)}^{2/\left(  3q-1\right)  }\Vert m\Vert
_{L^{2}(\Omega)}\leqslant\nonumber\\
&  \leqslant\frac{\varepsilon}{4}\Vert w^{q/2}\Vert_{W^{1,2}(\Omega)}%
^{2}+C_{15}\Vert m\Vert_{L^{2}(\Omega)}^{4\left(  3q-1\right)  }\Vert
w\Vert_{L^{1}(\Omega)}^{q}, \label{t2}%
\end{align}
where%
\[
C_{15}=C_{15}(\varepsilon,q,\Vert v_{0}\Vert_{L^{\infty}(\Omega)},\left\vert
\Omega\right\vert )=\left(  \frac{1}{3\left(  2q-1\right)  }\right)  \left(
\frac{\varepsilon}{4}\right)  ^{-3\left(  2q-1\right)  }\left(  \frac{3\left(
2q-1\right)  }{2(3q-1)}\left(  q-1\right)  C_{13}C(\Omega)\right)  ^{2\left(
3q-1\right)  }.
\]
From (\ref{es01}), using (\ref{mass4}), (\ref{mass6}), (\ref{t1}) and
(\ref{t2}), we get
\begin{align}
&  \frac{d}{dt}%
{\displaystyle\int\limits_{\Omega}}
z^{-1}w^{q}+\varepsilon\Vert w^{q/2}\Vert_{L^{2}(\Omega)}^{2}\leqslant
\nonumber\\
&  \leqslant-\frac{4(q-1)}{q}%
{\displaystyle\int\limits_{\Omega}}
z^{-1}\left\vert \nabla w^{q/2}\right\vert ^{2}+\frac{\varepsilon}{2}\Vert
w^{q/2}\Vert_{W^{1,2}(\Omega)}^{2}+C_{14}(\varepsilon)\Vert w\Vert
_{L^{1}(\Omega)}^{4q/\left(  3q-1\right)  }+C_{15}(\varepsilon)\Vert
m\Vert_{L^{2}(\Omega)}^{4\left(  3q-1\right)  }\Vert w\Vert_{L^{1}(\Omega
)}^{q}\leqslant\nonumber\\
&  \leqslant\left(  \frac{\varepsilon}{2}-\frac{4(q-1)}{q}\right)  \Vert\nabla
w^{q/2}\Vert_{L^{2}(\Omega)}^{2}+\frac{\varepsilon}{2}\Vert w^{q/2}%
\Vert_{L^{2}(\Omega)}^{2}+C_{16}, \label{es2}%
\end{align}
where%
\[
C_{16}=C_{16}(\varepsilon,q,\Vert v_{0}\Vert_{L^{\infty}(\Omega)},\left\vert
\Omega\right\vert )=\left(  C_{14}(\varepsilon,q)+C_{15}(\varepsilon
,q)C_{10}^{4\left(  3q-1\right)  }\right)  \left(  \max\{1,|\Omega|,\Vert
u_{0}\Vert_{L^{1}(\Omega)}\}\right)  ^{q}.
\]
Choosing $\varepsilon=2$, from (\ref{es2}) we obtain%
\[
\frac{d}{dt}%
{\displaystyle\int\limits_{\Omega}}
z^{-1}w^{q}+\frac{1}{C_{12}}%
{\displaystyle\int\limits_{\Omega}}
z^{-1}w^{q}\leqslant C_{16}.
\]
Applying Gronwall's lemma, the last inequality implies%
\[
\Vert w\Vert_{L^{q}(\Omega)}^{q}\leqslant%
{\displaystyle\int\limits_{\Omega}}
z^{-1}w^{q}\leqslant C_{12}\max\left\{  \Vert w_{0}\Vert_{L^{q}(\Omega)}%
^{q},C_{16}\right\}
\]
and we conclude the proof.
\end{proof}

\begin{corollary}
\label{cor_lp}Let $\Omega\subset\mathbb{R}^{3}$ be a domain with smooth
boundary and $(u_{0},v_{0},m_{0})\in L^{q}(\Omega)\times L^{\infty}%
(\Omega)\times W^{1,p}$, $1<q<\infty$, $p\geqslant\frac{6}{5}$. Then there
exists a constant $C_{17}=C_{17}(p,\left\Vert u_{0}\right\Vert _{L^{q}%
(\Omega)},\left\Vert v_{0}\right\Vert _{L^{\infty}(\Omega)},\left\Vert
m_{0}\right\Vert _{W^{1,p}(\Omega)})$ independent on time such that the
solution $u(x,t)$ to the system (\ref{e1})-(\ref{i}) satisfies \
\begin{equation}
\left\Vert u\right\Vert _{L^{\infty}(0,T;L^{q}(\Omega))}\leqslant
C_{17},\qquad\forall~1<q<+\infty. \label{lp_u}%
\end{equation}

\end{corollary}

We prove now the uniform boundedness for $w$ using the previous Proposition.
For this we use the iterative technique of Alikakos \cite{alikakos} and for
the sake of completeness this argument will be presented briefly in the next Proposition.

\begin{proposition}
\label{estimat_infinity}Let $\Omega\subset\mathbb{R}^{3}$ be a domain with
smooth boundary and $(w_{0},v_{0},m_{0})\in L^{\infty}(\Omega)\times
L^{\infty}(\Omega)\times W^{1,p}$, $p\geqslant\frac{6}{5}$. Then there exists
a constant $C_{18}=C_{18}(\left\Vert w_{0}\right\Vert _{L^{\infty}(\Omega
)},\left\Vert v_{0}\right\Vert _{L^{\infty}(\Omega)})$ independent on time
such that the solution $w(x,t)$ to the system (\ref{e1})-(\ref{i}) satisfies
\
\begin{equation}
\left\Vert w\right\Vert _{L^{\infty}(0,T;L^{\infty}(\Omega))}\leqslant C_{18}.
\label{l_inf}%
\end{equation}
\ \ \ \
\end{proposition}

\begin{proof}
We estimate the last term from the right-hand side of (\ref{es01}). Using
H\"{o}lder, Gagliardo-Nirenberg and Young's inequalities and (\ref{mass6})%
\begin{align}
\left(  q-1\right)  C_{13}%
{\displaystyle\int\limits_{\Omega}}
w^{q}m  &  \leqslant\left(  q-1\right)  C_{13}\Vert w^{q/2}\Vert_{L^{4}%
(\Omega)}^{2}\Vert m\Vert_{L^{2}(\Omega)}\leqslant\nonumber\\
&  \leqslant\left(  q-1\right)  C_{13}C(\Omega)\Vert w^{q/2}\Vert
_{W^{1,2}(\Omega)}^{9/5}\Vert w^{q/2}\Vert_{L^{1}(\Omega)}^{1/5}\Vert
m\Vert_{L^{2}(\Omega)}\leqslant\nonumber\\
&  \leqslant\varepsilon\Vert\nabla w^{q/2}\Vert_{L^{2}(\Omega)}^{2}%
+\varepsilon\Vert w^{q/2}\Vert_{L^{2}(\Omega)}^{2}+C_{19}(\varepsilon,q)\Vert
w^{q/2}\Vert_{L^{1}(\Omega)}^{2}, \label{es5}%
\end{align}
where%
\[
C_{19}=C_{19}(\varepsilon,q)=\frac{1}{9}\left(  \frac{9}{10}\right)
^{10}\varepsilon^{-9}\left[  \left(  q-1\right)  C_{13}C(\Omega)\right]
^{10}\Vert m\Vert_{L^{2}(\Omega)}^{10}.
\]
Introducing (\ref{es5}) in (\ref{es01}) we obtain%
\begin{align}
\frac{d}{dt}%
{\displaystyle\int\limits_{\Omega}}
z^{-1}w^{q}+\varepsilon\Vert w^{q/2}\Vert_{L^{2}(\Omega)}^{2}  &
\leqslant\left[  \varepsilon-\frac{4(q-1)}{p}\right]
{\displaystyle\int\limits_{\Omega}}
\left\vert \nabla w^{q/2}\right\vert ^{2}+\nonumber\\
&  +\left(  \mu qC_{12}+2\varepsilon\right)  \Vert w^{q/2}\Vert_{L^{2}%
(\Omega)}^{2}+C_{19}(\varepsilon,q)\Vert w^{q/2}\Vert_{L^{1}(\Omega)}^{2}.
\label{es6}%
\end{align}
Using again Gagliardo-Nirenberg and Young's inequalities, we infer%
\begin{align*}
\left\Vert w^{q/2}\right\Vert _{L^{2}(\Omega)}^{2}  &  \leqslant C_{20}%
(\Omega)\left\Vert w^{q/2}\right\Vert _{W^{1,2}(\Omega)}^{6/5}\left\Vert
w^{q/2}\right\Vert _{L^{1}(\Omega)}^{4/5}\leqslant\\
&  \leqslant\varepsilon\left\Vert \nabla w^{q/2}\right\Vert _{L^{2}(\Omega
)}^{2}+\varepsilon\left\Vert w^{q/2}\right\Vert _{L^{2}(\Omega)}^{2}%
+C_{21}\left(  \varepsilon\right)  \left\Vert w^{q/2}\right\Vert
_{L^{1}(\Omega)}^{2},
\end{align*}
or equivalently
\begin{equation}
\left\Vert w^{q/2}\right\Vert _{L^{2}(\Omega)}^{2}\leqslant\frac
{1}{1-\varepsilon}\left[  \varepsilon\left\Vert \nabla w^{q/2}\right\Vert
_{L^{2}(\Omega)}^{2}+C_{21}\left(  \varepsilon\right)  \left\Vert
w^{q/2}\right\Vert _{L^{1}(\Omega)}^{2}\right]  , \label{es7}%
\end{equation}
where $\varepsilon<1$ and%
\[
C_{21}=C_{21}\left(  \varepsilon\right)  =\frac{2}{3}\left(  \frac{3}%
{5}\right)  ^{5/2}\varepsilon^{-\frac{3}{2}}C_{20}^{5/2}(\Omega).
\]
Multiplying the inequality (\ref{es7}) by $\left(  \mu qC_{12}+2\varepsilon
\right)  $, we obtain from (\ref{es6})%
\begin{align}
\frac{d}{dt}%
{\displaystyle\int\limits_{\Omega}}
z^{-1}w^{q}+\varepsilon\Vert w^{q/2}\Vert_{L^{2}(\Omega)}^{2}  &
\leqslant\left[  \varepsilon-\frac{4(q-1)}{q}\right]
{\displaystyle\int\limits_{\Omega}}
\left\vert \nabla w^{q/2}\right\vert ^{2}+C_{19}(\varepsilon,q)\Vert
w^{q/2}\Vert_{L^{1}(\Omega)}^{2}+\nonumber\\
&  +\frac{\left(  \mu qC_{12}+2\varepsilon\right)  }{1-\varepsilon}\left[
\varepsilon\left\Vert \nabla w^{q/2}\right\Vert _{L^{2}(\Omega)}^{2}%
+C_{21}\left\Vert w^{q/2}\right\Vert _{L^{1}(\Omega)}^{2}\right]  =\nonumber\\
&  =\left[  \varepsilon+\frac{\left(  \mu qC_{12}+2\varepsilon\right)
}{1-\varepsilon}\varepsilon-\frac{4(q-1)}{q}\right]
{\displaystyle\int\limits_{\Omega}}
\left\vert \nabla w^{q/2}\right\vert ^{2}+\nonumber\\
&  +\left[  \frac{\left(  \mu qC_{12}+2\epsilon\right)  }{1-\varepsilon}%
C_{21}+C_{19}(\varepsilon,q)\right]  \left\Vert w^{q/2}\right\Vert
_{L^{1}(\Omega)}^{2}. \label{es8}%
\end{align}
We choose
\[
\varepsilon=\frac{1}{C_{22}\left(  q+1\right)  }<1
\]
where $C_{22}=\max\left\{  \mu C_{12},3\right\}  $. Hence
\[
\varepsilon+\frac{\left(  \mu qC_{12}+2\varepsilon\right)  }{1-\varepsilon
}\varepsilon<\frac{4(q-1)}{q}.
\]
Thus we obtain from (\ref{es8})%
\begin{equation}
\frac{d}{dt}%
{\displaystyle\int\limits_{\Omega}}
z^{-1}w^{q}+\frac{1}{C_{12}C_{22}\left(  q+1\right)  }%
{\displaystyle\int\limits_{\Omega}}
z^{-1}w^{q}\leqslant C_{23}(q)\left\Vert w^{q/2}\right\Vert _{L^{1}(\Omega
)}^{2}, \label{s9}%
\end{equation}
where%
\begin{align*}
C_{23}  &  =C_{23}(\varepsilon,q)=\frac{2}{3}\left(  \frac{3}{5}\right)
^{5/2}C_{20}^{5/2}(\Omega)\frac{\left(  \mu qC_{12}C_{22}\left(  q+1\right)
+2\right)  }{C_{22}\left(  q+1\right)  -1}\left(  C_{22}\left(  q+1\right)
\right)  ^{\frac{3}{2}}+\\
&  +\frac{1}{9}\left(  \frac{9}{10}\right)  ^{10}\left(  C_{22}\left(
q+1\right)  \right)  ^{9}\left[  \left(  q-1\right)  C_{13}C(\Omega)\right]
^{10}\Vert m\Vert_{L^{2}(\Omega)}^{10}.
\end{align*}
Applying Gronwall's lemma we deduce from the last inequality%
\begin{equation}%
{\displaystyle\int\limits_{\Omega}}
w^{q}\leqslant%
{\displaystyle\int\limits_{\Omega}}
z^{-1}w^{q}\leqslant C_{12}\max\left\{
{\displaystyle\int\limits_{\Omega}}
w_{0}^{q},C_{22}C_{23}(\varepsilon,q)\left(  q+1\right)  \left(
\underset{t\geqslant0}{\sup}%
{\displaystyle\int\limits_{\Omega}}
w^{q/2}\right)  ^{2}\right\}  . \label{es10}%
\end{equation}
We consider $q_{j}=2^{j}$, with $j\in\mathbb{N}$ and we use (\ref{es10})\ and
an iterative technique in order to estimate $\Vert w\Vert_{L^{q_{j}}(\Omega)}%
$. Denoting%
\[
a_{j}=C_{22}C_{23}(\varepsilon,q)\left(  q+1\right)  \leqslant C_{24}%
q_{j}^{20},
\]
the inequality (\ref{es10})\ becomes%
\begin{equation}%
{\displaystyle\int\limits_{\Omega}}
w^{q_{j}}\leqslant C_{12}\max\left\{  \max\left\{  \left\Vert w_{0}\right\Vert
_{L^{1}\left(  \Omega\right)  }^{q_{j}},\left\Vert w_{0}\right\Vert
_{L^{\infty}\left(  \Omega\right)  }^{q_{j}}\right\}  ,a_{j}\left(
\underset{t\geqslant0}{\sup}%
{\displaystyle\int\limits_{\Omega}}
w^{q_{j}/2}\right)  ^{2}\right\}  \label{es11}%
\end{equation}
for all $t>0$. Recursively, we obtain from (\ref{es11})%
\begin{align}
\left\Vert w\right\Vert _{L^{q_{j}}}  &  \leqslant C_{12}^{1/q_{j}}%
\max\left\{  \left\Vert w_{0}\right\Vert _{L^{1}\left(  \Omega\right)
},\left\Vert w_{0}\right\Vert _{L^{\infty}\left(  \Omega\right)  }\right\}
\left(  a_{j}a_{j-1}^{q_{1}}a_{j-2}^{q_{2}}...a_{1}^{q_{j-1}}\right)
^{\frac{1}{q_{j}}}\leqslant\nonumber\\
&  \leqslant2^{20\sum\limits_{k=1}^{j}\frac{j}{2^{j}}}C_{12}^{\frac{1}{2^{j}}%
}C_{24}^{\left(  1-\frac{1}{2^{j}}\right)  }\max\left\{  \left\Vert
w_{0}\right\Vert _{L^{1}\left(  \Omega\right)  },\left\Vert w_{0}\right\Vert
_{L^{\infty}\left(  \Omega\right)  }\right\}  . \label{es12}%
\end{align}
Taking $j\rightarrow\infty$ in the last inequality we deduce (\ref{l_inf}).
\end{proof}

\begin{corollary}
\label{cor_linf}Let $\Omega\subset\mathbb{R}^{3}$ be a domain with smooth
boundary and $(u_{0},v_{0},m_{0})\in L^{\infty}(\Omega)\times L^{\infty
}(\Omega)\times W^{1,p}(\Omega)$, $p\geqslant\frac{6}{5}$. Then there exists a
constant $C_{25}=C_{25}(\left\Vert u_{0}\right\Vert _{L^{\infty}(\Omega
)},\left\Vert v_{0}\right\Vert _{L^{\infty}(\Omega)})$ independent on time
such that the solution $u$ to the system (\ref{e1})-(\ref{i}) satisfies \
\begin{equation}
\left\Vert u\right\Vert _{L^{\infty}(0,T;L^{\infty}(\Omega))}\leqslant C_{25}.
\label{linf_u}%
\end{equation}

\end{corollary}

\begin{lemma}
\label{lemma_m}Let $\Omega\subset\mathbb{R}^{3}$ be a domain with smooth
boundary and $(u_{0},v_{0},m_{0})\in L^{p}(\Omega)\times L^{\infty}%
(\Omega)\times X_{p}^{\theta}(\Omega)$ where $p>3$, $\theta\in\left(  \frac
{1}{2},1\right)  $. Then%
\begin{equation}
\left\Vert m\left(  \cdot,t\right)  \right\Vert _{X_{p}^{\theta}}\leqslant
C_{26}, \label{m_est}%
\end{equation}
for all $t\in\left(  0,T\right)  $.
\end{lemma}

\begin{proof}
\noindent Taking into account the representation formula (\ref{rf}), the
hypothesis $(H_{2})$ and the estimate (\ref{lp_u}) we obtain
\begin{align*}
\left\Vert m\left(  \cdot,t\right)  \right\Vert _{X_{p}^{\theta}}  &
\leqslant C\Vert m_{0}\Vert_{X_{p}^{\theta}}+\int\limits_{0}^{t}\left\Vert
A_{2}^{\theta}e^{-\left(  t-s\right)  A_{2}}u(s)g(v(s))\right\Vert _{L^{p}%
}ds\leqslant\\
&  \leqslant C\Vert m_{0}\Vert_{X_{p}^{\theta}}+C(\theta)\int\limits_{0}%
^{t}(t-s)^{-\theta}e^{-\delta(t-s)}\left\Vert u(s)g(v(s))\right\Vert _{L^{p}%
}ds\leqslant\\
&  \leqslant C\Vert m_{0}\Vert_{X_{p}^{\theta}}+C_{17}C(\theta)\left(
L_{g}\left\Vert v_{0}\right\Vert _{L^{\infty}}+g(0)\right)  \int
\limits_{0}^{t}(t-s)^{-\theta}e^{-\delta(t-s)}ds.
\end{align*}
Denoting%
\[
C_{26}=\max\left\{  C\Vert m_{0}\Vert_{X_{p}^{\theta}},C_{17}C(\theta)\left(
L_{g}\left\Vert v_{0}\right\Vert _{L^{\infty}}+g(0)\right)  \delta^{\theta
-1}\Gamma(1-\theta)\right\}  ,
\]
where $\Gamma$ denotes the Gamma function, the last inequality implies
(\ref{m_est}).
\end{proof}

\begin{lemma}
\label{asdf}Let $\Omega\subset\mathbb{R}^{3}$ be a domain with smooth boundary
and $(u_{0},v_{0},m_{0})\in L^{p}(\Omega)\times W^{1,\infty}(\Omega)\times
X_{p}^{\theta}(\Omega)$ where $p>3$, $\theta\in\left(  \frac{3+p}%
{2p},1\right)  $. Then
\begin{equation}
\Vert v(t)\Vert_{W^{1,\infty}}\leqslant C_{27}\left(  1+t\right)  ,
\label{asdf1}%
\end{equation}
for all $t\in\left(  0,T\right)  $.
\end{lemma}

\begin{proof}
From the equation (\ref{e2}) we obtain
\begin{equation}
\nabla v=e^{-\int_{0}^{t}m}\left(  \nabla v_{0}-v_{0}\int_{0}^{t}\nabla
m\right)  \label{gh1}%
\end{equation}
which implies%
\[
\left\Vert \nabla v\right\Vert _{L^{\infty}}\leqslant\left\Vert \nabla
v_{0}\right\Vert _{L^{\infty}}+\left\Vert v_{0}\right\Vert _{L^{\infty}}%
\int_{0}^{t}\left\Vert m\right\Vert _{W^{1,\infty}}ds\leqslant\left\Vert
v_{0}\right\Vert _{w^{1,\infty}}\max\left\{  1,C_{26}\right\}  \left(
1+t\right)  .
\]
Next, taking into account Lemma \ref{lemma_m} and denoting%
\[
C_{27}=\Vert v_{0}\Vert_{L^{\infty}(\Omega)}\max\left\{  1,C_{26}\right\}
\]
we conclude the proof.
\end{proof}

\begin{lemma}
\label{lemma_u}Let $\Omega\subset\mathbb{R}^{3}$ be a domain with smooth
boundary and $(u_{0},v_{0},m_{0})\in W^{1,q}(\Omega)\times W^{1,\infty}%
(\Omega)\times X_{p}^{\theta}(\Omega)$ where $p>3$, $\theta\in\left(
\frac{3+p}{2p},1\right)  $ and $q\geqslant\frac{3p}{p+3}$. For \ $t\in(0,T)$
\ we have
\begin{equation}
\Vert u(\cdot,t)\Vert_{W^{1,q}(\Omega)}\leqslant C_{28}(1+t)e^{C_{29}t},
\label{cot2}%
\end{equation}
for all $t\in\left(  0,T\right)  $.
\end{lemma}

\begin{proof}
First we establish a bound for $\Vert w(\cdot,t)\Vert_{W^{1,q}(\Omega)}$.
Taking into account the representation formula%
\begin{equation}
w(x,t)=e^{-tA_{1}}w_{0}+\int\limits_{0}^{t}e^{-\left(  t-s\right)  A_{1}}%
G_{1}(w,v,m)(s)ds,\quad\forall t\in(0,T) \label{rm}%
\end{equation}
we get%
\begin{align*}
&  \Vert w(\cdot,t)\Vert_{W^{1,q}(\Omega)}\leqslant\Vert e^{-tA_{1}}w_{0}%
\Vert_{W^{1,q}(\Omega)}+\int\limits_{0}^{t}\left\Vert A_{1}^{\theta
}e^{-\left(  t-s\right)  A_{1}}G_{1}(w,v,m)\right\Vert _{L^{q}\left(
\Omega\right)  }ds\leqslant\\
&  \leqslant C\left\Vert w_{0}\right\Vert _{W^{1,q}(\Omega)}+C\left(
\theta,q\right)  C_{27}\left[  L_{\chi}\left\Vert v_{0}\right\Vert
_{L^{\infty}\left(  \Omega\right)  }+\chi(0)\right]  \int\limits_{0}%
^{t}\left(  t-s\right)  ^{-\theta}e^{-\delta(t-s)}\left(  1+s\right)
\left\Vert w\right\Vert _{W^{1,q}\left(  \Omega\right)  }ds+\\
&  +C_{11}\left[  \left(  \mu+1\right)  +\mu C_{11}C_{12}+\mu\left\Vert
v_{0}\right\Vert _{L^{\infty}\left(  \Omega\right)  }\right]  \int
\limits_{0}^{t}\left(  t-s\right)  ^{-\theta}e^{-\delta(t-s)}ds+\\
&  +C_{11}C_{26}L_{\chi}\left\Vert v_{0}\right\Vert _{L^{\infty}\left(
\Omega\right)  }\left[  \left\Vert v_{0}\right\Vert _{L^{\infty}\left(
\Omega\right)  }+\chi(0)\right]  \int\limits_{0}^{t}\left(  t-s\right)
^{-\theta}e^{-\delta(t-s)}ds\leqslant\\
&  \leqslant C_{30}+C_{31}\int\limits_{0}^{t}\left(  t-s\right)  ^{-\theta
}e^{-\delta(t-s)}\left(  1+s\right)  \left\Vert w\right\Vert _{W^{1,q}\left(
\Omega\right)  }ds
\end{align*}
where%
\begin{align*}
C_{30}  &  =C\left\Vert w_{0}\right\Vert _{W^{1,q}(\Omega)}+C_{11}%
\delta^{\theta-1}\max\left\{  C_{26}\left\Vert v_{0}\right\Vert _{L^{\infty
}\left(  \Omega\right)  }\left[  L_{\chi}\left\Vert v_{0}\right\Vert
_{L^{\infty}\left(  \Omega\right)  }+\chi(0)\right]  \right.  ,\\
&  \left.  \left[  \left(  \mu+1\right)  +\mu C_{11}C_{12}+\mu\left\Vert
v_{0}\right\Vert _{L^{\infty}\left(  \Omega\right)  }\right]  \right\}
\Gamma(1-\theta),\\
C_{31}  &  =C\left(  \theta,q\right)  C_{27}\left[  L_{\chi}\left\Vert
v_{0}\right\Vert _{L^{\infty}\left(  \Omega\right)  }+\chi(0)\right]  .
\end{align*}
and $\Gamma$ is the Gamma function. Applying Gronwall's lemma in the
last
inequality we obtain%
\begin{equation}
\Vert w(\cdot,t)\Vert_{W^{1,q}(\Omega)}\leqslant C_{30}e^{\Gamma\left(
1-\theta\right)  C_{31}}e^{\Gamma\left(  1-\theta\right)  C_{31}t}%
.\noindent\label{wp}%
\end{equation}
Finally, as $u=wz^{-1}$ we get from (\ref{asdf1}) and (\ref{wp})
\begin{align*}
\Vert u(\cdot,t)\Vert_{W^{1,q}(\Omega)}  &  \leqslant\Vert z^{-1}%
\Vert_{W^{1,\infty}(\Omega)}\Vert w\Vert_{W^{1,q}(\Omega)}\leqslant\\
&  \leqslant C_{12}\left[  1+C_{27}\left(  L_{\chi}\left\Vert v_{0}\right\Vert
_{L^{\infty}\left(  \Omega\right)  }+\chi(0)\right)  (1+t)\right]  \Vert
w\Vert_{W^{1,q}(\Omega)}\leqslant\\
&  \leqslant2C_{12}C_{30}\max\left\{  1,C_{27}\left(  L_{\chi}\left\Vert
v_{0}\right\Vert _{L^{\infty}\left(  \Omega\right)  }+\chi(0)\right)
\right\}  e^{\Gamma\left(  1-\theta\right)  C_{31}}\left(  1+t\right)
e^{\Gamma\left(  1-\theta\right)  C_{31}t}.
\end{align*}
Denoting
\begin{align*}
C_{28}  &  =2C_{12}C_{30}\max\left\{  1,C_{27}\left(  L_{\chi}\left\Vert
v_{0}\right\Vert _{L^{\infty}\left(  \Omega\right)  }+\chi(0)\right)
\right\}  e^{C_{31}\Gamma\left(  1-\theta\right)  },\\
C_{29}  &  =C_{31}\Gamma\left(  1-\theta\right)
\end{align*}
the last inequality implies (\ref{cot2}).
\end{proof}

\section{Steady-states}

\label{steadyst}

\setcounter{equation}{0} \setcounter{figure}{0}

Additionally, we suppose that the function $g$ satisfies the following
property:%
\[
(H_{3})\qquad g(v)\neq0\text{ if }v\neq0\text{.}%
\]

In this section we deal with the stationary problem associated to
(\ref{e1})-(\ref{i}). More precisely, we are interested in the solution of the
system%
\begin{align}
&  0=\Delta w+\chi(v)\nabla v\cdot\nabla w+\mu w\left(  1-wz^{-1}-v\right)  &
x  &  \in\Omega\label{s1}\\
&  0=mv & x  &  \in\Omega\label{s2}\\
&  0=d\Delta m-\gamma m+wz^{-1}g(v) & x  &  \in\Omega\label{s3}\\
&  \frac{\partial w}{\partial\eta}=\frac{\partial m}{\partial\eta}=0 & x  &
\in\partial\Omega. \label{s4}%
\end{align}

\begin{theorem}
\label{steadys}If $(w^{\ast},v^{\ast},m^{\ast})\in C^{1}(\overline{\Omega
})\times W^{1,\infty}(\Omega)\times C^{1}(\overline{\Omega})$ are the
non-negative solutions to (\ref{s1})-(\ref{s4}), then they are given by
\begin{align*}
(w^{\ast},v^{\ast},m^{\ast})  &  =(0,\widetilde{v},0)\,,\\
(w^{\ast},v^{\ast},m^{\ast})  &  =(k,0,k\gamma^{-1}g(0))\,
\end{align*}
where $k=0$ or $k=1$ if $\mu>0$, $k\geqslant0$ is an arbitrary constant if
$\mu=0$ and $\widetilde{v}\in W^{1,\infty}(\Omega)$ is an arbitrary
non-negative function.
\end{theorem}

\begin{proof}
We distinguish between two cases: $w\left(  x\right)  g(v\left(  x\right)
)>0$ at some point $x\in\Omega$, or $wg(v)\equiv0$.

First, we suppose that there exists at least a point $x\in\Omega$ such that
$w\left(  x\right)  g(v\left(  x\right)  )\neq0$. We claim that
\[
\min_{x\in\overline{\Omega}}m(x)>0.
\]
To prove this, we assume that there exists $x_{0}\in\overline{\Omega}$ such
that $m(x_{0})=0$. We have two possibilities:

\begin{enumerate}
\item[a.] if $x_{0}\in\Omega$, then \cite[Theorem 3.5]{gilbarg} implies that
$m$ is a constant function and\ we deduce that $m(x)\equiv0$ for all
$x\in\Omega$, but this is not a solution for the equation (\ref{s3}).

\item[b.] if $x_{0}\in\partial\Omega$, then using \cite[Lemma 3.4]{gilbarg} we
get
\[
\frac{\partial m}{\partial\eta}(x_{0})<0,
\]
which contradicts (\ref{s4}).
\end{enumerate}

Since $\min\limits_{x\in\overline{\Omega}}m(x)>0$, the equation (\ref{s2})
implies $v(x)\equiv0$ for all $x\in\Omega$. Therefore the function $w(x)$ is a
solution of the following equation
\begin{equation}
\left\{
\begin{array}
[c]{ll}%
-\Delta w=\mu w(1-w) & x\in\Omega\\
\displaystyle\frac{\partial w}{\partial\eta}=0 & x\in\partial\Omega.
\end{array}
\right.  \label{st-1}%
\end{equation}

\begin{enumerate}
\item[a.] If $\mu=0$, we deduce that the solutions of (\ref{st-1}) are
$w\equiv k$ where $k$ is a non-negative constant.

\item[b.] If $\mu>0$, then using a similar argument as in \cite{BO86} we
obtain that $w\equiv0$ and $w\equiv1$ are the only non-negative solutions to
(\ref{st-1}).
\end{enumerate}

Therefore, from (\ref{s3})-(\ref{s4}) we obtain that $m$ is a constant and
moreover
\[
m\equiv k\gamma^{-1}g(0)
\]
where $k$ is a non-negative constant if $\mu=0$ and $k=0$ or $k=1$ if $\mu>0$.

We suppose now that for all $x\in\Omega$, $wg(v)\equiv0$. Then from
(\ref{s3})-(\ref{s4}) we get $m\equiv0$.

We observe that $wg(v)\equiv0$ implies $wv\equiv0,$ using the hypothesis
$(H_{3})$. Arguing as in the previous case we prove that either $w\equiv0$ or
$\min\limits_{x\in\overline{\Omega}}w(x)>0$.

\begin{enumerate}
\item[a.] If $w\equiv0$ then any $\widetilde{v}\in W^{1,\infty}(\Omega)$ ,
$\widetilde{v}\geqslant0$ solves (\ref{s2}).

\item[b.] If $\min\limits_{x\in\overline{\Omega}}w(x)>0$ then $wg(v)\equiv0$
implies $g(v)\equiv0$ and by $(H_{3})$ $v\equiv0$. Therefore, (\ref{s1}) with
the boundary condition can be written in the form (\ref{st-1}). Hence $w\equiv
k$ where $k$ is a non-negative constant if $\mu=0$ and $k=0$ or $k=1$ if
$\mu>0$.
\end{enumerate}
\end{proof}

\section{Asymptotic behaviour of global solutions}

\setcounter{equation}{0} \setcounter{figure}{0}

In what follows we study the asymptotic behaviour of the global
solution $(u,v,m)$ to the problem (\ref{eq1})-(\ref{ic}). Under some
additionally hypotheses on the initial data we shall prove the
convergence of the solutions to the steady states.

Hereafter $\Omega\subset\mathbb{R}^{3}$ is a domain with smooth boundary.

\begin{lemma}
\label{bound_1}Let $\mu\geqslant0$, $v_{0}\in L^{\infty}(\Omega)$ be positive
and if $\mu>0$ we assume $0<v_{0}(x)<1$ for all $x\in\Omega$. If there exists
a constant $a>0$ such that $u_{0}(x)\geqslant a$, then every global solution
$u(x,t)$ satisfies
\begin{equation}
u(x,t)\geqslant\min\left\{  1,a\right\}  e^{\int\limits_{0}^{v_{0}}\chi(s)ds},
\label{u_bound}%
\end{equation}
for all $x\in\Omega$, $t>0$. Moreover, if there exists a positive constant $M$
such that
\begin{equation}
g(v)\geqslant M \label{hg}%
\end{equation}
for all $v\in\mathbb{R}_{+}$, then there exists a constant $\sigma>0$ such
that
\begin{equation}
m(x,t)\geqslant\sigma>0, \label{m_delta}%
\end{equation}
for all $x\in\Omega$, $t>0$.
\end{lemma}

\begin{proof}
Let $\rho$ be a positive constant to be chosen later. By multiplying the
equation (\ref{eq1}) with $z^{-1}(w-\rho)_{-}$, where $(w-\rho)_{-}%
=\max\left\{  \rho-w,0\right\}  $, and after that integrating over $\Omega$ we
get%
\begin{align}
\frac{1}{2}\frac{d}{dt}\int\limits_{\Omega}\left[  z^{-1}(w-\rho)_{-}%
^{2}\right]   &  =-\int\limits_{\Omega}z^{-1}\left\vert \nabla(w-\rho
)_{-}\right\vert ^{2}-\frac{1}{2}\int\limits_{\Omega}z^{-1}\chi(v)mv(w-\rho
)_{-}^{2}-\nonumber\\
&  -\int\limits_{\Omega}z^{-1}\chi(v)mvw(w-\rho)_{-}-\mu\int\limits_{\Omega
}z^{-1}w\left(  1-wz^{-1}-v\right)  (w-\rho)_{-}. \label{sigma}%
\end{align}
Let us notice that in the case $\mu=0$ the right-hand side of (\ref{sigma}) is
non-positive. If $\mu>0$, choosing the constant $0<\rho\leqslant\left(
1-v_{0}\right)  e^{-\int\limits_{0}^{v_{0}}\chi(s)ds}$ implies that the last
term in (\ref{sigma}) is also non-positive. Integrating (\ref{sigma}) in time
and taking into account the above considerations we have
\[
\int\limits_{\Omega}z^{-1}(w-\rho)_{-}^{2}\leqslant\int\limits_{\Omega}%
e^{\int\limits_{0}^{v_{0}}\chi(s)ds}(w_{0}-\rho)_{-}^{2}.
\]
Using the same procedure as in \cite{limo}\ finally we obtain
\begin{equation}
u(x,t)\geqslant\min\left\{  1,a\right\}  e^{\int\limits_{0}^{v_{0}}\chi(s)ds},
\label{sigma1}%
\end{equation}
for all $x\in\Omega$, $t>0$. Next, we conclude the proof by applying the
maximum principle for the parabolic problem%
\begin{align*}
&  m_{t}-d\Delta m+\gamma m=ug(v) & x  &  \in\Omega,\quad t\in\mathbb{R}_{+}\\
&  \frac{\partial m}{\partial\eta}=0 & x  &  \in\partial\Omega,\quad
t\in\mathbb{R}_{+}\\
&  m(x,0)=m_{0}(x) & x  &  \in\Omega
\end{align*}
using the hypothesis $g(v)\geqslant M>0$ and the estimate (\ref{sigma1}).
\end{proof}

\begin{lemma}
\label{bound_2}Let $\left(  u_{0},v_{0},m_{0}\right)  \in L^{1}(\Omega
)\times\left(  W^{1,2}(\Omega)\cap L^{\infty}(\Omega)\right)  \times
W^{1,p}(\Omega)$, $p\geqslant\frac{6}{5}$. Assume that $v_{0}>0$ and the
assumption (\ref{hg}) is satisfied. Then we have
\begin{equation}
\int_{\Omega}|\nabla\left(  v^{q/2}(t)\right)  |^{2}\leqslant C_{33}e^{-kt},
\label{dv_min}%
\end{equation}
for all $t>0$, $q\geqslant1$ and $0<k<q\sigma$.
\end{lemma}

\begin{proof}
Using the equation (\ref{eq2}) we deduce%
\begin{equation}
\frac{d}{dt}\int_{\Omega}|\nabla v^{q/2}|^{2}=-q\int_{\Omega}m\left\vert
\nabla v^{q/2}\right\vert ^{2}-\frac{q}{2}\int_{\Omega}\nabla\left(
v^{q}\right)  \cdot\nabla m, \label{est30}%
\end{equation}
for $q\geqslant1$. On the other hand, multiplying (\ref{eq3}) with $v^{q}$,
$q\geqslant1$ and integrating in space, we have%
\begin{equation}
-\frac{d}{2}\int_{\Omega}\nabla m\cdot\nabla\left(  v^{q}\right)  =\frac{1}%
{2}\frac{d}{dt}\int_{\Omega}mv^{q}+\frac{q}{2}\int_{\Omega}m^{2}v^{q}%
+\frac{\gamma}{2}\int_{\Omega}mv^{q}-\frac{1}{2}\int_{\Omega}ug(v)v^{q}.
\label{est32}%
\end{equation}
Inserting (\ref{est32}) in (\ref{est30}) we get
\begin{equation}
\frac{d}{dt}\int_{\Omega}|\nabla v^{q/2}|^{2}=-q\int_{\Omega}m\left\vert
\nabla v^{q/2}\right\vert ^{2}+\frac{q}{2d}\frac{d}{dt}\int_{\Omega}%
mv^{q}+\frac{q^{2}}{2d}\int_{\Omega}m^{2}v^{q}+\frac{\gamma q}{2d}\int
_{\Omega}mv^{q}-\frac{q}{2d}\int_{\Omega}ug(v)v^{q}. \label{est34}%
\end{equation}
Now, taking into account (\ref{m_delta}) and multiplying the equality
(\ref{est34}) by $e^{kt}$, $0<k<q\sigma$, we obtain%
\begin{equation}
\frac{d}{dt}\left(  e^{kt}\int_{\Omega}|\nabla v^{q/2}|^{2}\right)
\leqslant\frac{q}{2d}\frac{d}{dt}\left(  e^{kt}\int\limits_{\Omega}%
mv^{q}\right)  +\frac{q^{2}}{2d}e^{kt}\int\limits_{\Omega}m^{2}v^{q}%
+\frac{\gamma q}{2d}e^{kt}\int\limits_{\Omega}mv^{q}. \label{est36}%
\end{equation}
Integrating the last inequality on $(0,t)$ we have%
\begin{equation}
e^{kt}\int\limits_{\Omega}|\nabla v^{q/2}|^{2}\leqslant\int\limits_{\Omega
}|\nabla v_{0}^{q/2}|^{2}+\frac{q}{2d}e^{kt}\int\limits_{\Omega}mv^{q}%
+\frac{q^{2}}{2d}\int\limits_{0}^{t}e^{ks}\int\limits_{\Omega}m^{2}v^{q}%
+\frac{\gamma q}{2d}\int\limits_{0}^{t}e^{ks}\int\limits_{\Omega}mv^{q}.
\label{est37}%
\end{equation}
Since $v(x,t)=v_{0}e^{-\int_{0}^{t}m}$ and $m(x,t)\geqslant\sigma>0$,
\begin{equation}
v(x,t)\leqslant\left\Vert v_{0}\right\Vert _{L^{\infty}\left(  \Omega\right)
}e^{-\sigma t}. \label{vest}%
\end{equation}
Introducing the estimate (\ref{vest}) in (\ref{est37}) and taking into account
(\ref{mass3}) and (\ref{mass6}), we get
\begin{align}
&  \int\limits_{\Omega}|\nabla v^{q/2}(t)|^{2}\leqslant e^{-kt}\int
\limits_{\Omega}|\nabla v_{0}^{q/2}|^{2}+\nonumber\\
&  +\frac{q}{2d}\left\Vert v_{0}\right\Vert _{L^{\infty}\left(  \Omega\right)
}^{q}\max\left\{  \left\Vert m\right\Vert _{L^{1}(\Omega)},C_{10}\right\}
\left(  e^{\left(  k-\sigma q\right)  t}+q\int\limits_{0}^{t}e^{\left(
k-\sigma q\right)  s}+\gamma\int\limits_{0}^{t}e^{\left(  k-\sigma q\right)
s}\right)  e^{-kt}\leqslant\nonumber\\
&  \leqslant\left[  \int\limits_{\Omega}|\nabla v_{0}^{q/2}|^{2}+\frac{q}%
{2d}C_{32}\left\Vert v_{0}\right\Vert _{L^{\infty}\left(  \Omega\right)  }%
^{q}\left(  1+\frac{q+\gamma}{\sigma q-k}\right)  \right]  e^{-kt},
\label{dvest}%
\end{align}
where%
\[
C_{32}=\max\left\{  C_{10},\Vert m_{0}\Vert_{L^{1}(\Omega)}+\left(  L_{g}\Vert
v_{0}\Vert_{L^{\infty}(\Omega)}+g(0)\right)  \max\{|\Omega|,\Vert u_{0}%
\Vert_{L^{1}(\Omega)}\}\right\}  .
\]
Finally, from (\ref{dvest}) and denoting%
\[
C_{33}=\left[  \int\limits_{\Omega}|\nabla v_{0}^{q/2}|^{2}+\frac{q}{2d}%
C_{32}\left\Vert v_{0}\right\Vert _{L^{\infty}\left(  \Omega\right)  }%
^{q}\left(  1+\frac{q+\gamma}{\sigma q-k}\right)  \right]  ,
\]
we conclude the proof.
\end{proof}

\begin{theorem}
\label{bound_3}If the hypotheses of Lemma \ref{bound_2} are satisfied, then
\begin{align}
&  \Vert u(\cdot,t)-\overline{u}\Vert_{L^{p}(\Omega)}\leqslant C_{36}%
e^{-C_{37}t},\label{asymp1}\\
&  \Vert v(\cdot,t)\Vert_{L^{\infty}(\Omega)}\leqslant C_{38}e^{-\sigma
t},\label{asymp2}\\
&  \Vert m(\cdot,t)-\overline{u}\gamma^{-1}g(0)\Vert_{L^{p}(\Omega)}\leqslant
C_{39}e^{-C_{40}t}, \label{asymp3}%
\end{align}
where $C_{36}$, $C_{37}$, $C_{39}$, $C_{40}$ are positive constants
independent on $t$, $C_{38}=\Vert v_{0}\Vert_{L^{\infty}(\Omega)}$ and%
\begin{equation}
\overline{u}=\left\{
\begin{array}
[c]{ll}%
{\frac{1}{|\Omega|}}\int\limits_{\Omega}u_{0}, & \mu=0\\
1, & \mu>0.
\end{array}
\right.  \label{u_int}%
\end{equation}

\end{theorem}

\begin{proof}
Let $\alpha$ be a positive constant to be chosen later. By multiplying the
equation (\ref{eq1}) with $\left(  u-\alpha\right)  ^{2p+1}$, $p\geqslant0$
and integrating in space we get%
\begin{align}
&  \frac{d}{dt}\int\limits_{\Omega}(u-\alpha)^{2p+2}\leqslant-\left(
p+1\right)  \left(  2p+1\right)  \int\limits_{\Omega}\left(  u-\alpha\right)
^{2p}\left\vert \nabla u\right\vert ^{2}+\nonumber\\
&  +2^{2p}\left(  p+1\right)  \left(  2p+1\right)  C_{25}\max\left\{
C_{25},\alpha^{2p}\right\}  \left(  L_{\chi}\left\Vert v_{0}\right\Vert
_{L^{\infty}(\Omega)}^{2}+\chi(0)\right)  ^{2}\int\limits_{\Omega}\left\vert
\nabla v\right\vert ^{2}+\nonumber\\
&  +2\mu\left(  p+1\right)  \int\limits_{\Omega}u\left(  1-u\right)  \left(
u-\alpha\right)  ^{2p+1}-2\mu\left(  p+1\right)  \int\limits_{\Omega}%
u^{2}v\left(  u-\alpha\right)  ^{2p}+2\mu\alpha\left(  p+1\right)
\int\limits_{\Omega}uv\left(  u-\alpha\right)  ^{2p}. \label{est1}%
\end{align}
In the case $\mu=0$ we consider $\alpha=\overline{u}={\frac{1}{|\Omega|}}%
\int\limits_{\Omega}u_{0}$. Using the Poincar\'{e} inequality we obtain from
(\ref{est1})%
\[
\frac{d}{dt}\int\limits_{\Omega}(u-\overline{u})^{2}+C_{34}\int\limits_{\Omega
}(u-\overline{u})^{2}\leqslant\left\Vert u\right\Vert _{L^{\infty
}(0,t;L^{\infty}(\Omega))}^{2}\left(  L_{\chi}\left\Vert v_{0}\right\Vert
_{L^{\infty}(\Omega)}^{2}+\chi(0)\right)  ^{2}\int\limits_{\Omega}\left\vert
\nabla v\right\vert ^{2}.
\]
Using Lemma \ref{bound_2} for $q=1$, $k\neq C_{34}$ and applying the Gronwall
inequality in the last estimate we have
\begin{equation}
\Vert u(t)-\overline{u}\Vert_{L^{2}(\Omega)}^{2}\leqslant\left[  \Vert
u_{0}-\overline{u}\Vert_{L^{2}(\Omega)}^{2}+C_{25}^{2}\left(  L_{\chi
}\left\Vert v_{0}\right\Vert _{L^{\infty}(\Omega)}^{2}+\chi(0)\right)
^{2}\frac{C_{33}}{\left\vert k-C_{34}\right\vert }\right]  e^{-\min\left\{
k,C_{34}\right\}  t}. \label{u_bar}%
\end{equation}
Moreover, for $p\geqslant2$ it follows that
\begin{align}
\Vert u(t)-\overline{u}\Vert_{L^{p}(\Omega)}^{p}  &  \leqslant2^{p-2}%
\max\left\{  \overline{u}^{p-2},C_{25}^{p-2}\right\}  \left[  \Vert
u_{0}-\overline{u}\Vert_{L^{2}(\Omega)}^{2}\right.  +\nonumber\\
&  +\left.  C_{25}^{2}\left(  L_{\chi}\left\Vert v_{0}\right\Vert _{L^{\infty
}(\Omega)}^{2}+\chi(0)\right)  ^{2}\frac{C_{33}}{\left\vert k-C_{34}%
\right\vert }\right]  e^{-\min\left\{  k,C_{34}\right\}  t}. \label{lpu1}%
\end{align}
In the case $\mu>0$ we consider $\alpha=1$ and taking into account Lemma
\ref{bound_1} and Lemma \ref{bound_2}, we obtain from (\ref{est1})
\begin{align*}
&  \frac{d}{dt}\int\limits_{\Omega}(u-1)^{2p+2}+C_{35}\int\limits_{\Omega
}\left(  u-1\right)  ^{2p+2}\leqslant\\
&  \leqslant2^{2p}C_{25}\max\left\{  1,C_{25}^{2p}\right\}  \left[  \left(
p+1\right)  \left(  2p+1\right)  C_{25}\left(  L_{\chi}\left\Vert
v_{0}\right\Vert _{L^{\infty}(\Omega)}^{2}+\chi(0)\right)  ^{2}\int
\limits_{\Omega}\left\vert \nabla v\right\vert ^{2}+2\mu\left\vert
\Omega\right\vert \left\Vert v_{0}\right\Vert _{L^{\infty}(\Omega)}e^{-\delta
t}\right]
\end{align*}
where
\[
C_{35}=2\mu\left(  p+1\right)  \min\left\{  1,a\right\}  e^{\int
\limits_{0}^{v_{0}}\chi(s)ds}.
\]
Applying Gronwall's inequality, using the estimate (\ref{dv_min}) and choosing
$k\neq C_{35}$ we obtain\
\begin{align}
&  \Vert u(t)-1\Vert_{L^{2p+2}(\Omega)}^{2p+2}\leqslant\left\{  \Vert
u_{0}-1\Vert_{L^{2p+2}(\Omega)}^{2p+2}+2^{2p}C_{25}\max\left\{  1,C_{25}%
^{2p}\right\}  \left[  2\mu\left\vert \Omega\right\vert \left\Vert
v_{0}\right\Vert _{L^{\infty}(\Omega)}+\right.  \right. \nonumber\\
&  +\left.  \left.  \left(  p+1\right)  \left(  2p+1\right)  C_{25}\left(
L_{\chi}\left\Vert v_{0}\right\Vert _{L^{\infty}(\Omega)}^{2}+\chi(0)\right)
^{2}\frac{C_{33}}{\left\vert k-C_{35}\right\vert }\right]  \right\}
e^{-\min\left\{  k,C_{35}\right\}  t}. \label{uone}%
\end{align}
Taking into account (\ref{lpu1}) and (\ref{uone}) we conclude the estimate
(\ref{asymp1}).

By multiplying the equation (\ref{eq1}) with $\left(  m-\overline{u}%
\gamma^{-1}g(0)\right)  ^{2p+1}$, $p\geqslant0$ and integrating in space we
get
\begin{align*}
&  \frac{d}{dt}\int\limits_{\Omega}(m-\overline{u}\gamma^{-1}g(0))^{2p+2}%
=-2d\left(  p+1\right)  ^{2}\int\limits_{\Omega}(m-\overline{u}\gamma
^{-1}g(0))^{2p}\left\vert \nabla m\right\vert ^{2}-\\
&  -2\gamma\left(  p+1\right)  \int\limits_{\Omega}(m-\overline{u}\gamma
^{-1}g(0))^{2p+2}+2\left(  p+1\right)  \int\limits_{\Omega}(m-\overline
{u}\gamma^{-1}g(0))^{2p+1}\left(  ug(v)-\overline{u}g(0)\right)  \leqslant\\
&  \leqslant-2d\left(  p+1\right)  ^{2}\int\limits_{\Omega}(m-\overline
{u}\gamma^{-1}g(0))^{2p}\left\vert \nabla m\right\vert ^{2}-\gamma\left(
p+1\right)  \int\limits_{\Omega}(m-\overline{u}\gamma^{-1}g(0))^{2p+2}+\\
&  +\frac{\left(  2p+1\right)  ^{2p+1}}{\gamma^{2p+1}\left(  p+1\right)
^{4p+3}}\int\limits_{\Omega}\max\left\{  L_{g}u^{2p+2}v^{2p+2},g^{2p+2}%
(0)\left\vert u-\overline{u}\right\vert ^{2p+2}\right\}  \leqslant\\
&  \leqslant-\gamma\left(  p+1\right)  \int\limits_{\Omega}(m-\overline
{u}\gamma^{-1}g(0))^{2p+2}+\\
&  +\frac{\left(  2p+1\right)  ^{2p+1}}{\gamma^{2p+1}\left(  p+1\right)
^{4p+3}}\max\left\{  L_{g}\left\vert \Omega\right\vert \left(  C_{17}%
\left\Vert v(\cdot,t)\right\Vert _{L^{\infty}(\Omega)}\right)  ^{2p+2}%
,g^{2p+2}(0)\left\Vert u\left(  \cdot,t\right)  -\overline{u}\right\Vert
_{L^{2p+2}}^{2p+2}\right\}  .
\end{align*}
Using the estimates (\ref{asymp1}), (\ref{asymp2}) and the Gronwall
lemma, when
$\gamma\neq\min\left\{  C_{37},2\sigma\right\}  $ we obtain%
\begin{align*}
&  \left\Vert m-\overline{u}\gamma^{-1}g(0)\right\Vert _{L^{2p+2}}%
^{2p+2}\leqslant\left[  \left\Vert m_{0}-\overline{u}\gamma^{-1}%
g(0)\right\Vert _{L^{2p+2}}^{2p+2}\right.  +\\
&  +\left.  \frac{\left(  2p+1\right)  ^{2p+1}\max\left\{  L_{g}\left\vert
\Omega\right\vert \left(  C_{17}C_{38}\right)  ^{2p+2},g^{2p+2}(0)C_{36}%
\right\}  }{\gamma^{2p+1}\left(  p+1\right)  ^{4p+4}\left\vert \gamma
-\min\left\{  C_{37},2\sigma\right\}  \right\vert }\right]  e^{-\min\left\{
C_{37},\left(  2p+2\right)  \sigma,\gamma(p+1)\right\}  t}%
\end{align*}
and we conclude the estimate (\ref{asymp3}).
\end{proof}

In the next theorem we shall prove that starting with initial data in suitable
spaces we obtain stronger convergences.

\begin{theorem}
\label{bound_4}Let $\left(  u_{0},v_{0},m_{0}\right)  \in L^{\infty}%
(\Omega)\times W^{1,\infty}(\Omega)\times W^{1,p}(\Omega)$, $p>3$. Under the
hypotheses of Lemma \ref{bound_1}, we have%
\begin{align}
&  \Vert v(\cdot,t)\Vert_{W^{1,\infty}}\leqslant C_{41}\left(  1+t\right)
e^{-\sigma t},\label{asy2}\\
&  \Vert m(\cdot,t)-\overline{u}\gamma^{-1}g(0)\Vert_{X_{p}^{\theta}}\leqslant
C_{42}\left(  1+t\right)  t^{-\theta}e^{-C_{43}t}, \label{asy3}%
\end{align}
for all $t\in\left(  0,t\right)  $. Moreover, if the hypotheses of Lemma
\ref{lemma_u} are satisfied, then%
\begin{equation}
\Vert u(\cdot,t)-\overline{u}\Vert_{W^{1,\infty}}\leqslant C_{44}e^{-\frac
{1}{2}C_{43}t}, \label{asy1}%
\end{equation}
for $t$ sufficiently large. The constants $C_{41}$, $C_{42}$, $C_{43}$,
$C_{44}$ are positive and independent on $t$, and $\overline{u}$ is given by
(\ref{u_int}).
\end{theorem}

\begin{proof}
Let $\theta\in\left(  \frac{3+p}{2p},1\right)  $. From the representation
formula (\ref{rf}) and taking into account (\ref{h3}), (\ref{lp_u}),
(\ref{asymp1}), (\ref{asymp2}) we obtain
\begin{align*}
&  \left\Vert m\left(  \cdot,t\right)  -\overline{u}\gamma^{-1}g(0)\right\Vert
_{X_{p}^{\theta}}\leqslant\Vert e^{-tA_{2}}\left[  m_{0}-\overline{u}%
\gamma^{-1}g(0)\right]  \Vert_{X_{p}^{\theta}}+\int\limits_{0}^{t}\left\Vert
A_{2}^{\theta}e^{-\left(  t-s\right)  A_{2}}\left[  ug(v)-\overline
{u}g(0)\right]  \right\Vert _{L^{p}\left(  \Omega\right)  }ds\leqslant\\
&  \leqslant C(\theta)t^{-\theta}e^{-\delta t}\left(  \Vert m_{0}\Vert
_{L^{p}(\Omega)}+\Vert\overline{u}\gamma^{-1}g(0)\Vert_{L^{p}(\Omega)}\right)
+\\
&  +C\left(  \theta,p\right)  \int\limits_{0}^{t}\left(  t-s\right)
^{-\theta}e^{-\delta\left(  t-s\right)  }\left[  \left\Vert
ug(v)-ug(0)\right\Vert _{L^{p}\left(  \Omega\right)  }+\left\Vert
ug(0)-\overline{u}g(0)\right\Vert _{L^{p}\left(  \Omega\right)  }\right]
ds\leqslant\\
&  \leqslant C(\theta)t^{-\theta}e^{-\delta t}\left(  \Vert m_{0}\Vert
_{L^{p}(\Omega)}+\overline{u}\gamma^{-1}g(0)\left\vert \Omega\right\vert
\right)  +\\
&  +C\left(  \theta,p\right)  L_{g}C_{17}C_{38}\int\limits_{0}^{t}\left(
t-s\right)  ^{-\theta}e^{-\delta\left(  t-s\right)  }e^{-\sigma s}ds+C\left(
\theta,p\right)  C_{36}g(0)\int\limits_{0}^{t}\left(  t-s\right)  ^{-\theta
}e^{-\delta\left(  t-s\right)  }e^{-C_{37}s}ds\leqslant\\
&  \leqslant\max\left\{  C(\theta)\left(  \Vert m_{0}\Vert_{L^{p}(\Omega
)}+\overline{u}\gamma^{-1}g(0)\left\vert \Omega\right\vert \right)
,\frac{C\left(  \theta,p\right)  }{1-\theta}\left(  L_{g}C_{17}C_{38}%
+g(0)C_{36}\right)  \right\}  \left(  1+t\right)  t^{-\theta}e^{-\min\left\{
\sigma,\delta,C_{37}\right\}  t}.
\end{align*}
Denoting
\begin{align*}
C_{42}  &  =\max\left\{  C(\theta)\left(  \Vert m_{0}\Vert_{L^{p}(\Omega
)}+\overline{u}\gamma^{-1}g(0)\left\vert \Omega\right\vert \right)
,\frac{C\left(  \theta,p\right)  }{1-\theta}\left(  L_{g}C_{17}C_{38}%
+g(0)C_{36}\right)  \right\}  ,\\
C_{43}  &  =\min\left\{  \sigma,\delta,C_{37}\right\}
\end{align*}
the last inequality implies the estimate (\ref{asy3}).

Taking into account (\ref{gh1}), the estimate (\ref{asy1}) and the embedding
(\ref{h2}), it follows that%
\begin{align*}
\left\Vert \nabla v\right\Vert _{L^{\infty}}  &  \leqslant e^{-\sigma
t}\left(  \left\Vert \nabla v_{0}\right\Vert _{L^{\infty}}+\left\Vert
v_{0}\right\Vert _{L^{\infty}}\int_{0}^{t}\left\Vert m-\overline{u}\gamma
^{-1}g(0)\right\Vert _{W^{1,\infty}}ds\right)  \leqslant\\
&  \leqslant e^{-\sigma t}\left(  \left\Vert \nabla v_{0}\right\Vert
_{L^{\infty}}+C_{42}\left\Vert v_{0}\right\Vert _{L^{\infty}}\int_{0}%
^{t}(1+s)s^{-\theta}e^{-C_{43}s}ds\right)  \leqslant\\
&  \leqslant2\max\left\{  1,C_{42}C_{43}^{\theta-1}\Gamma\left(
1-\theta\right)  \right\}  \left\Vert v_{0}\right\Vert _{W^{1,\infty}%
}(1+t)e^{-\sigma t}.
\end{align*}
The last inequality together with (\ref{asymp2}) imply (\ref{asy2}) where%
\[
C_{41}=\left(  1+2\max\left\{  1,C_{42}C_{43}^{\theta-1}\Gamma\left(
1-\theta\right)  \right\}  \right)  \left\Vert v_{0}\right\Vert _{W^{1,\infty
}}.
\]
We start now to prove (\ref{asy1}). Let $t_{0}>0$ sufficiently large. Taking
into account the embedding (\ref{h2}), we obtain from (\ref{asy3})%
\[
\left\Vert \nabla m\right\Vert _{L^{\infty}\left(  \Omega\right)  }\leqslant
C_{45}e^{-\frac{1}{2}C_{43}t}%
\]
for all $t>t_{0}$. If we take also into account (\ref{gh1}), the last
inequality implies%
\begin{align}
\left\Vert \nabla v\right\Vert _{L^{\infty}\left(  \Omega\right)  }  &
\leqslant e^{-\sigma t}\left(  \left\Vert \nabla v_{0}\right\Vert _{L^{\infty
}\left(  \Omega\right)  }+C_{45}\left\Vert v_{0}\right\Vert _{L^{\infty
}\left(  \Omega\right)  }\int_{0}^{t}e^{-\frac{1}{2}C_{43}t}\right)
\leqslant\nonumber\\
&  \leqslant\max\left\{  1,\frac{2C_{45}}{C_{43}}\right\}  \left\Vert
v_{0}\right\Vert _{W^{1,\infty}\left(  \Omega\right)  }e^{-\sigma t}
\label{dv}%
\end{align}
for $t$ sufficiently large. From the representation formula%
\[
w(x,t)=e^{-\left(  t-t_{0}\right)  A_{1}}w(t_{0})+\int\limits_{t_{0}}%
^{t}e^{-\left(  t-s\right)  A_{1}}G_{1}(w,v,m)(s)ds
\]
and taking into account (\ref{lp}), (\ref{l_inf}) and (\ref{asymp2}), we
obtain \ \
\begin{align}
&  \left\Vert w\left(  \cdot,t\right)  -\overline{u}\right\Vert _{X_{p}%
^{\theta}}\leqslant\Vert e^{-\left(  t-t_{0}\right)  A_{1}}\left(
w(t_{0})-\overline{u}\right)  \Vert_{X_{p}^{\theta}}+\nonumber\\
&  +\int\limits_{t_{0}}^{t}\left\Vert A_{1}^{\theta}e^{-\left(  t-s\right)
A_{1}}\left[  \chi(v)\nabla v\cdot\nabla w+\left(  \mu+1\right)  w-\mu
w\left(  wz^{-1}+v\right)  +\chi(v)wvm-\overline{u}\right]  \right\Vert
_{L^{p}\left(  \Omega\right)  }ds\leqslant\nonumber\\
&  \leqslant C(\theta)\left(  t-t_{0}\right)  ^{-\theta}e^{-\delta\left(
t-t_{0}\right)  }\Vert w(t_{0})-\overline{u}\Vert_{L^{p}(\Omega)}+\nonumber\\
&  +C\left(  \theta,p\right)  \left(  L_{\chi}\left\Vert v_{0}\right\Vert
_{L^{\infty}\left(  \Omega\right)  }+\chi(0)\right)  \int\limits_{t_{0}}%
^{t}\left(  t-s\right)  ^{-\theta}e^{-\delta\left(  t-s\right)  }\left\Vert
\nabla v\right\Vert _{L^{\infty}\left(  \Omega\right)  }\left\Vert \nabla
w\right\Vert _{L^{p}\left(  \Omega\right)  }ds+\nonumber\\
&  +\mu C_{18}C\left(  \theta,p\right)  \int\limits_{t_{0}}^{t}\left(
t-s\right)  ^{-\theta}e^{-\delta\left(  t-s\right)  }\left\Vert u-1\right\Vert
_{L^{p}\left(  \Omega\right)  }+C\left(  \theta,p\right)  \int\limits_{t_{0}%
}^{t}\left(  t-s\right)  ^{-\theta}e^{-\delta\left(  t-s\right)  }\left\Vert
w-\overline{u}\right\Vert _{L^{p}\left(  \Omega\right)  }+\nonumber\\
&  +\mu C_{11}C_{38}C\left(  \theta,p\right)  \int\limits_{t_{0}}^{t}\left(
t-s\right)  ^{-\theta}e^{-\delta\left(  t-s\right)  }e^{-\sigma s}%
ds+\nonumber\\
&  +C_{11}C_{38}C\left(  \theta,p\right)  \left(  L_{\chi}\left\Vert
v_{0}\right\Vert _{L^{\infty}\left(  \Omega\right)  }+\chi(0)\right)
\int\limits_{t_{0}}^{t}\left(  t-s\right)  ^{-\theta}e^{-\delta\left(
t-s\right)  }e^{-\sigma s}\left\Vert m\right\Vert _{L^{\infty}\left(
\Omega\right)  }. \label{w1}%
\end{align}
Let us observe that%
\[
\left\Vert w-u\right\Vert _{L^{p}\left(  \Omega\right)  }^{p}\leqslant
\left\Vert u\right\Vert _{L^{\infty}\left(  \Omega\right)  }^{p}\left\Vert
z-1\right\Vert _{L^{p}\left(  \Omega\right)  }^{p}\leqslant C_{46}\left\Vert
u\right\Vert _{L^{\infty}\left(  \Omega\right)  }^{p}\int\limits_{\Omega
}\left\vert \int_{0}^{v}\chi(s)ds\right\vert ^{p}\leqslant C_{47}%
^{p}e^{-\sigma pt},
\]
where%
\[
C_{47}=2C_{46}^{1/p}C_{25}C_{38}\left\vert \Omega\right\vert ^{1/p}%
\max\left\{  C_{38}\left(  \frac{L_{\chi}}{2}\right)  ,\chi(0)\right\}  .
\]
From the above inequality, using (\ref{asymp1}), it follows also that
\begin{equation}
\left\Vert w-\overline{u}\right\Vert _{L^{p}\left(  \Omega\right)  }%
\leqslant\max\left\{  C_{36},C_{47}\right\}  e^{-\min\left\{  \sigma
,C_{37}\right\}  t}. \label{wu}%
\end{equation}
Using the estimates (\ref{m_est}), (\ref{asymp1}), (\ref{dv}), (\ref{asymp2})
and (\ref{wu}), we obtain from (\ref{w1})
\begin{align*}
\left\Vert w\left(  \cdot,t\right)  -\overline{u}\right\Vert _{X_{p}^{\theta
}}  &  \leqslant C_{48}\left(  t-t_{0}\right)  ^{1-\theta}e^{-\min\left\{
\delta,\sigma,C_{37}\right\}  t}+C_{49}\left(  t-t_{0}\right)  ^{-\theta
}e^{-\delta t}+\\
&  +C_{50}\int\limits_{t_{0}}^{t}\left(  t-s\right)  ^{-\theta}e^{-\delta
\left(  t-s\right)  }e^{-\sigma s}\left\Vert \nabla w\right\Vert
_{L^{p}\left(  \Omega\right)  }ds
\end{align*}
where%
\begin{align*}
C_{48}  &  =\frac{C\left(  \theta,p\right)  }{1-\theta}\left\{  \left[
\max\left\{  C_{36},C_{47}\right\}  +C_{11}C_{26}C_{38}\left(  L_{\chi
}\left\Vert v_{0}\right\Vert _{L^{\infty}\left(  \Omega\right)  }%
+\chi(0)\right)  \right]  +\mu\left(  C_{18}C_{36}+C_{11}C_{38}\right)
\right\} \\
C_{49}  &  =C(\theta)e^{\delta t_{0}}\left(  C_{11}+\Vert\overline{u}%
\Vert_{L^{p}(\Omega)}\right) \\
C_{50}  &  =C\left(  \theta,p\right)  \left(  L_{\chi}\left\Vert
v_{0}\right\Vert _{L^{\infty}\left(  \Omega\right)  }+\chi(0)\right)
\max\left\{  1,\frac{2C_{42}}{C_{43}}\right\}  \left\Vert v_{0}\right\Vert
_{W^{1,\infty}\left(  \Omega\right)  }.
\end{align*}
Taking into account the embedding (\ref{h1}) and applying the Gronwall lemma
we obtain from the last inequality%
\begin{align}
&  \left\Vert w\left(  \cdot,t\right)  -\overline{u}\right\Vert _{X_{p}%
^{\theta}}\leqslant C_{48}\left(  t-t_{0}\right)  ^{1-\theta}e^{-\min\left\{
\delta,\sigma,C_{37}\right\}  t}+C_{49}\left(  t-t_{0}\right)  ^{-\theta
}e^{-\delta t}+\nonumber\\
&  +C_{50}e^{C_{50}\delta^{\theta-1}\Gamma(1-\theta)}\left[  C_{48}%
e^{-\min\left\{  \delta,\sigma,C_{37}\right\}  t}\int\limits_{t_{0}}%
^{t}\left(  s-t_{0}\right)  ^{1-\theta}\left(  t-s\right)  ^{-\theta}%
ds+C_{49}e^{-\delta t}\int\limits_{t_{0}}^{t}\left(  s-t_{0}\right)
^{-\theta}\left(  t-s\right)  ^{-\theta}ds\right]  \leqslant\nonumber\\
&  \leqslant C_{45}\left(  t-t_{0}\right)  ^{1-\theta}e^{-\min\left\{
\delta,\sigma,C_{37}\right\}  t}+C_{49}\left(  t-t_{0}\right)  ^{-\theta
}e^{-\delta t}+\nonumber\\
&  +C_{50}e^{C_{50}\delta^{\theta-1}\Gamma(1-\theta)}\frac{2^{2\theta-2}%
\sqrt{\pi}\left(  t-t_{0}\right)  ^{1-2\theta}\Gamma(1-\theta)}{\Gamma
(\frac{3}{2}-\theta)}\left(  C_{48}\left(  t-t_{0}\right)  +2C_{49}\right)
e^{-\min\left\{  \delta,\sigma,C_{37}\right\}  t}\leqslant\nonumber\\
&  \leqslant\left(  t-t_{0}\right)  ^{-\theta}\left(  C_{48}\left(
t-t_{0}\right)  +2C_{49}\right)  \left(  1+2^{2\theta-2}C_{50}\frac{\sqrt{\pi
}\left(  t-t_{0}\right)  ^{1-\theta}\Gamma(1-\theta)}{\Gamma(\frac{3}%
{2}-\theta)}e^{C_{50}\delta^{\theta-1}\Gamma(1-\theta)}\right)  e^{-\min
\left\{  \delta,\sigma,C_{37}\right\}  t}. \label{w13}%
\end{align}
Let us observe that%
\begin{equation}
\Vert u(\cdot,t)-\overline{u}\Vert_{W^{1,\infty}}\leqslant\Vert z^{-1}%
\Vert_{W^{1,\infty}}\cdot\Vert w-\overline{u}\Vert_{W^{1,\infty}}%
+\Vert\overline{u}\Vert_{L^{\infty}}\Vert z^{-1}-1\Vert_{W^{1,\infty}}.
\label{ubar}%
\end{equation}
Using (\ref{dv}) we obtain%
\begin{equation}
\Vert z^{-1}\Vert_{W^{1,\infty}}\leqslant C_{12}\left[  1+C_{12}\left(
L_{\chi}\left\Vert v_{0}\right\Vert _{L^{\infty}\left(  \Omega\right)  }%
+\chi(0)\right)  \max\left\{  1,\frac{2C_{45}}{C_{43}}\right\}  \left\Vert
v_{0}\right\Vert _{W^{1,\infty}\left(  \Omega\right)  }\right]  e^{-\sigma t}.
\label{star1}%
\end{equation}
Taking into account (\ref{asymp2}) and (\ref{dv}), we estimate
\begin{align}
&  \Vert z^{-1}-1\Vert_{W^{1,\infty}}\leqslant C_{46}\Vert\int_{0}^{v}%
\chi(s)ds\Vert_{L^{\infty}}+C_{12}\left(  L_{\chi}\left\Vert v_{0}\right\Vert
_{L^{\infty}\left(  \Omega\right)  }+\chi(0)\right)  \max\left\{
1,\frac{2C_{45}}{C_{43}}\right\}  \left\Vert v_{0}\right\Vert _{W^{1,\infty
}\left(  \Omega\right)  }e^{-\sigma t}\leqslant\nonumber\\
&  \leqslant\left[  C_{46}C_{38}\left(  \frac{L_{\chi}}{2}C_{38}%
+\chi(0)\right)  +C_{12}\left(  L_{\chi}\left\Vert v_{0}\right\Vert
_{L^{\infty}\left(  \Omega\right)  }+\chi(0)\right)  \max\left\{
1,\frac{2C_{45}}{C_{43}}\right\}  \left\Vert v_{0}\right\Vert _{W^{1,\infty
}\left(  \Omega\right)  }\right]  e^{-\sigma t} \label{star2}%
\end{align}
From (\ref{ubar}), using the embedding (\ref{h1}) and the estimates
(\ref{w13}), (\ref{star1}), (\ref{star2}), we obtain
\begin{align*}
&  \Vert u(t)-\overline{u}\Vert_{W^{1,\infty}}\leqslant C_{51}\Vert
w-\overline{u}\Vert_{W^{1,\infty}}+C_{52}\overline{u}e^{-\sigma t}\leqslant\\
&  \leqslant\left[  C_{51}\left(  t-t_{0}\right)  ^{-\theta}\left(
C_{48}\left(  t-t_{0}\right)  +2C_{49}\right)  \left(  1+2^{2\theta-2}%
C_{50}\frac{\sqrt{\pi}\left(  t-t_{0}\right)  ^{1-\theta}\Gamma(1-\theta
)}{\Gamma(\frac{3}{2}-\theta)}e^{C_{50}\delta^{\theta-1}\Gamma(1-\theta
)}\right)  +C_{52}\overline{u}\right]  e^{-\min\left\{  \delta,\sigma
,C_{37}\right\}  t}%
\end{align*}
where
\begin{align*}
C_{51}  &  =C_{12}\left[  1+C_{12}\left(  L_{\chi}\left\Vert v_{0}\right\Vert
_{L^{\infty}\left(  \Omega\right)  }+\chi(0)\right)  \max\left\{
1,\frac{2C_{45}}{C_{43}}\right\}  \left\Vert v_{0}\right\Vert _{W^{1,\infty
}\left(  \Omega\right)  }\right] \\
C_{52}  &  =\left[  C_{46}C_{38}\left(  \frac{L_{\chi}}{2}C_{38}%
+\chi(0)\right)  +C_{12}\left(  L_{\chi}\left\Vert v_{0}\right\Vert
_{L^{\infty}\left(  \Omega\right)  }+\chi(0)\right)  \max\left\{
1,\frac{2C_{45}}{C_{43}}\right\}  \left\Vert v_{0}\right\Vert _{W^{1,\infty
}\left(  \Omega\right)  }\right]
\end{align*}
The last inequality implies (\ref{asy1}) for $t$ sufficiently large.
\end{proof}

In what follows we investigate the case when $g(0)=0$ and we start with an
auxiliary lemma.

\begin{lemma}
\label{lemmaf}Let $f\in C^{1}(0,+\infty)$ satisfying%
\[
\int\limits_{0}^{\infty}\left\vert f(s)\right\vert ds\leqslant C,\quad
\int\limits_{0}^{\infty}\left\vert f^{\prime}(s)\right\vert ds\leqslant C
\]
then%
\[
\lim_{t\rightarrow\infty}f(t)=0.
\]

\end{lemma}

\begin{proof}
Assume that $\lim\limits_{t\rightarrow\infty}f(t)\neq0$, then there exists a
sequence $\{t_{n}\}_{n\in\mathbb{N}}$, $t_{n}\rightarrow\infty$ such that%
\[
\left\vert f(t_{n})\right\vert >C>0,\quad\forall n\geqslant n_{0}.
\]
We can assume that $t_{n+1}>t_{n}+1$ (otherwise we take a subsequence). Let
$0<\tau<1$, then for all $n\geqslant n_{0}$%
\[
\left\vert \left\vert f(t_{n}+\tau)\right\vert -\left\vert f(t_{n})\right\vert
\right\vert \leqslant\left\vert f(t_{n}+\tau)-f(t_{n})\right\vert =\left\vert
\int\limits_{t_{n}}^{t_{n}+\tau}f^{\prime}(s)ds\right\vert \leqslant
\int\limits_{t_{n}}^{t_{n}+\tau}\left\vert f^{\prime}(s)\right\vert
ds\leqslant\int\limits_{t_{n}}^{t_{n}+1}\left\vert f^{\prime}(s)\right\vert
ds.
\]
But
\[
\lim_{n\rightarrow\infty}\int\limits_{t_{n}}^{t_{n}+1}\left\vert f^{\prime
}(s)\right\vert ds=0
\]
(otherwise, on a subsequence, $\int\limits_{t_{n}}^{t_{n}+1}\left\vert
f^{\prime}(s)\right\vert ds\geqslant C$\ which implies $\int\limits_{0}%
^{\infty}\left\vert f^{\prime}(s)\right\vert ds\geqslant\sum\limits_{n}%
\int\limits_{t_{n}}^{t_{n}+1}\left\vert f^{\prime}(s)\right\vert
ds\rightarrow\infty$). Then for $n\geqslant n_{1}\geqslant n_{0}$ we have%
\[
\int\limits_{t_{n}}^{t_{n}+1}\left\vert f^{\prime}(s)\right\vert ds<\frac
{C}{2}.
\]
Therefore
\[
\left\vert f(t_{n}+\tau)\right\vert >\frac{C}{2},
\]
which implies \
\[
\int\limits_{t_{n}}^{t_{n}+1}\left\vert f(s)\right\vert ds\geqslant\frac{C}%
{2},
\]
which contradicts the hypothesis.
\end{proof}

\begin{theorem}
\label{bound_5}Let $\mu>0$ and $(u_{0},v_{0},m_{0})$ $\in L^{1}(\Omega
)\times\left(  W^{1,2}(\Omega)\cap L^{\infty}(\Omega)\right)  \times
W^{1,p}(\Omega)$, $p\geqslant\frac{6}{5}$. Let $g$ be a function satisfying
the hypotheses $(H2)$\ and $(H3)$and\ $g(0)=0$. Then
\begin{align}
&  \lim_{t\rightarrow\infty}\Vert m(\cdot,t)\Vert_{L^{2}(\Omega)}%
=0,\label{asymp6}\\
&  \lim_{t\rightarrow\infty}\Vert g(v(\cdot,t))v^{q}(\cdot,t)\Vert
_{L^{1}(\Omega)}=0,\qquad q\geqslant1. \label{asymp5}%
\end{align}
Moreover, if $g(v)\geqslant C_{53}v$ for all $v\geqslant0$, where $C_{53}$ is
a positive constant independent on $t$, then%
\begin{equation}
\lim_{t\rightarrow\infty}\Vert u(\cdot,t)-\overline{u}\Vert_{L^{2}(\Omega)}=0,
\label{asymp4}%
\end{equation}
where $\overline{u}$ is given by (\ref{u_int}).
\end{theorem}

\begin{proof}
On multiplying the equation (\ref{eq3}) by $m$, integrating over $\Omega$ and
taking into account the estimates (\ref{v}), (\ref{l_inf}) and the hypothesis
$(H_{3})$, we obtain%
\begin{equation}
\frac{1}{2}\frac{d}{dt}\int\limits_{\Omega}m^{2}+d\int\limits_{\Omega
}\left\vert \nabla m\right\vert ^{2}+\gamma\int\limits_{\Omega}m^{2}%
=-L_{g}\int\limits_{\Omega}uv_{t}. \label{m2}%
\end{equation}
Integrating in time the last equality we get%
\[
\frac{1}{2}\int\limits_{\Omega}m^{2}+d\int\limits_{0}^{t}\int\limits_{\Omega
}\left\vert \nabla m\right\vert ^{2}+\gamma\int\limits_{0}^{t}\int
\limits_{\Omega}m^{2}\leqslant C_{25}L_{g}\int\limits_{\Omega}v_{0}+\frac
{1}{2}\int\limits_{\Omega}m_{0}^{2}.
\]
Therefore we have%
\begin{align}
\int\limits_{0}^{t}\int\limits_{\Omega}\left\vert \nabla m\right\vert ^{2}  &
\leqslant C,\label{dm}\\
\int\limits_{0}^{t}\int\limits_{\Omega}m^{2}  &  \leqslant C. \label{m_estim}%
\end{align}
From (\ref{m2}) we deduce%
\[
\left\vert \frac{d}{dt}\int\limits_{\Omega}m^{2}\right\vert \leqslant
-2C_{12}C_{18}L_{g}\int\limits_{\Omega}v_{t}+2d\int\limits_{\Omega}\left\vert
\nabla m\right\vert ^{2}+2\gamma\int\limits_{\Omega}m^{2}%
\]
and integrating in time we obtain%
\begin{equation}
\int\limits_{0}^{t}\left\vert \frac{d}{dt}\int\limits_{\Omega}m^{2}\right\vert
\leqslant-2C_{12}C_{18}L_{g}\int\limits_{0}^{t}\int\limits_{\Omega}%
v_{t}+2d\int\limits_{0}^{t}\int\limits_{\Omega}\left\vert \nabla m\right\vert
^{2}+2\gamma\int\limits_{0}^{t}\int\limits_{\Omega}m^{2}. \label{dm_estim}%
\end{equation}
The last inequality together with (\ref{dm}), (\ref{m_estim}) and Lemma
\ref{lemmaf} imply (\ref{asymp6}).

In order to deal with the convergence of $v$, we integrate in time the
estimate (\ref{est34}) and we obtain%
\begin{align*}
&  \int_{\Omega}|\nabla v^{q/2}|^{2}+q\int\limits_{0}^{t}\int_{\Omega}m|\nabla
v^{q/2}|^{2}+\frac{q}{2d}\int\limits_{0}^{t}\int_{\Omega}ug(v)v^{q}\leqslant\\
&  \leqslant\int_{\Omega}|\nabla v_{0}^{q/2}|^{2}+\frac{q}{2d}\int_{\Omega
}mv^{q}+\frac{q^{2}}{2d}\int\limits_{0}^{t}\int_{\Omega}m^{2}v^{q}%
+\frac{\gamma q}{2d}\int\limits_{0}^{t}\int_{\Omega}mv^{q}.
\end{align*}
The last inequality, taking into account also (\ref{u_bound}), implies%
\begin{equation}
\int\limits_{0}^{t}\int_{\Omega}g(v)v^{q}\leqslant C_{54}. \label{gv1}%
\end{equation}
On the other hand%
\begin{align}
\int\limits_{0}^{t}\left\vert \frac{d}{dt}\int\limits_{\Omega}g(v)v^{q}%
\right\vert  &  \leqslant-L_{g^{\prime}}\int\limits_{0}^{t}\int\limits_{\Omega
}v^{q+1}v_{t}-\left\vert g^{\prime}(0)\right\vert \int\limits_{0}^{t}%
\int\limits_{\Omega}v^{q}v_{t}-qL_{g}\int\limits_{0}^{t}\int\limits_{\Omega
}v^{q}v_{t}\leqslant\nonumber\\
&  \leqslant\left(  L_{g^{\prime}}\left\Vert v_{0}\right\Vert _{L^{\infty
}(\Omega)}+\left\vert g^{\prime}(0)\right\vert +qL_{g}\right)  \left\vert
\Omega\right\vert \left\Vert v_{0}\right\Vert _{L^{\infty}(\Omega)}%
^{q+1}\leqslant C. \label{gv2}%
\end{align}
From Lemma \ref{lemmaf}, (\ref{gv1}) and (\ref{gv2}) we get (\ref{asymp5}).

In order to have the last estimate we multiply the equation (\ref{e1}) by
$(w-a)$, where $a$ is a constant to be determined later and we obtain%
\begin{align}
\frac{1}{2}\frac{d}{dt}\int\limits_{\Omega}\left[  z^{-1}(w-a)^{2}\right]   &
=-\int\limits_{\Omega}z^{-1}\left\vert \nabla w\right\vert ^{2}+\mu
\int\limits_{\Omega}wz^{-1}(1-w)(w-a)+\mu\int\limits_{\Omega}w^{2}%
z^{-1}(1-z^{-1})(w-a)-\nonumber\\
&  -\mu\int\limits_{\Omega}wz^{-1}v(w-a)+\int\limits_{\Omega}\chi
(v)w^{2}z^{-1}mv-a\int\limits_{\Omega}\chi(v)wz^{-1}mv. \label{z}%
\end{align}
We consider first the case $\mu=0$. Taking $a=\overline{u}=\frac{1}{\left\vert
\Omega\right\vert }\int\limits_{\Omega}u$, from the equation (\ref{e2}) and
integrating in time the last inequality we get
\begin{align*}
&  \int\limits_{\Omega}z^{-1}(w-\overline{u})^{2}+2\int\limits_{0}^{t}%
\int\limits_{\Omega}z^{-1}\left\vert \nabla w\right\vert ^{2}\leqslant
\int\limits_{\Omega}\left[  z_{0}^{-1}(w_{0}-\overline{u})^{2}\right]
+2\int\limits_{0}^{t}\int\limits_{\Omega}w^{2}z^{-1}\chi(v)mv\leqslant\\
&  \leqslant\int\limits_{\Omega}\left[  z_{0}^{-1}(w_{0}-\overline{u}%
)^{2}\right]  +2C_{12}\left\vert \Omega\right\vert \left(  L_{\chi}\left\Vert
v_{0}\right\Vert _{L^{\infty}(\Omega)}+\chi(0)\right)  \left\Vert
v_{0}\right\Vert _{L^{\infty}(\Omega)}\left\Vert w\right\Vert _{L^{\infty
}(0,t;L^{\infty}(\Omega))}^{2}.
\end{align*}
Therefore we have
\begin{equation}
\int\limits_{0}^{t}\int\limits_{\Omega}\left\vert \nabla w\right\vert
^{2}\leqslant\int\limits_{0}^{t}\int\limits_{\Omega}z^{-1}\left\vert \nabla
w\right\vert ^{2}\leqslant C_{55}. \label{dw1}%
\end{equation}
From (\ref{z}) we deduce%
\[
\left\vert \frac{d}{dt}\int\limits_{\Omega}\left[  z^{-1}(w-\overline{u}%
)^{2}\right]  \right\vert \leqslant2\int\limits_{\Omega}z^{-1}\left\vert
\nabla w\right\vert ^{2}-2\int\limits_{\Omega}\chi(v)w^{2}z^{-1}%
v_{t}-2\overline{u}\int\limits_{\Omega}\chi(v)wz^{-1}v_{t}.
\]
Integrating in time the last inequality and taking into account (\ref{mass1}),
(\ref{l_inf}) and (\ref{dw1}) we obtain
\begin{equation}
\int\limits_{0}^{t}\left\vert \frac{d}{dt}\int\limits_{\Omega}\left[
z^{-1}(w-\overline{u})^{2}\right]  \right\vert \leqslant2\int\limits_{0}%
^{t}\int\limits_{\Omega}z^{-1}\left\vert \nabla w\right\vert ^{2}%
-2\int\limits_{0}^{t}\int\limits_{\Omega}\chi(v)w^{2}z^{-1}v_{t}-2\overline
{u}\int\limits_{0}^{t}\int\limits_{\Omega}\chi(v)wz^{-1}v_{t}\leqslant C_{56}
\label{d1}%
\end{equation}
where
\[
C_{56}=2C_{55}+2C_{12}C_{18}\left\vert \Omega\right\vert ^{2}\left(
\max\left\{  1,\Vert u_{0}\Vert_{L^{\infty}(\Omega)}\right\}  +C_{18}\right)
\left(  L_{\chi}\left\Vert v_{0}\right\Vert _{L^{\infty}\left(  \Omega\right)
}+\chi(0)\right)  \left\Vert v_{0}\right\Vert _{L^{\infty}\left(
\Omega\right)  }.
\]
Let us remark that we have the following estimate by using Poincar\'{e}
inequality and (\ref{dw1})%
\begin{align}
\int\limits_{0}^{t}\int\limits_{\Omega}z^{-1}(w-\overline{u})^{2}  &
\leqslant2C_{12}\int\limits_{0}^{t}\int\limits_{\Omega}(w-\overline{w}%
)^{2}+2C_{12}\int\limits_{0}^{t}\int\limits_{\Omega}(\overline{w}-\overline
{u})^{2}\leqslant\nonumber\\
&  \leqslant2CC_{12}\int\limits_{0}^{t}\int\limits_{\Omega}\left\vert \nabla
w\right\vert ^{2}+2C_{12}\int\limits_{0}^{t}\int\limits_{\Omega}(\overline
{w}-\overline{u})^{2}\leqslant C+2C_{12}\int\limits_{0}^{t}\int\limits_{\Omega
}(\overline{w}-\overline{u})^{2} \label{d2}%
\end{align}
where $\overline{w}=\frac{1}{\left\vert \Omega\right\vert }\int\limits_{\Omega
}w$. Then, in order to obtain the estimate (\ref{asymp4}) using Lemma
\ref{lemmaf}, it is enough to prove that%
\begin{equation}
\int\limits_{0}^{t}\int\limits_{\Omega}(\overline{w}-\overline{u}%
)^{2}\leqslant C. \label{wu_bar}%
\end{equation}
Using (\ref{gv1}) for $q=1$, an easy calculation shows us that%
\begin{equation}
\int\limits_{0}^{t}\int\limits_{\Omega}(1-z^{-1})^{2}\leqslant C_{12}%
\max_{s\in(0,\left\Vert v_{0}\right\Vert _{L^{\infty}\left(  \Omega\right)
})}\chi^{2}(s)\int\limits_{0}^{t}\int\limits_{\Omega}v^{2}\leqslant
C_{12}C_{53}C_{54}\max_{s\in(0,\left\Vert v_{0}\right\Vert _{L^{\infty}\left(
\Omega\right)  })}\chi^{2}(s). \label{sh}%
\end{equation}
Taking into account the last estimate we obtain%
\begin{align}
\int\limits_{0}^{t}\int\limits_{\Omega}(\overline{w}-\overline{u})^{2}  &
\leqslant\int\limits_{0}^{t}\int\limits_{\Omega}\frac{1}{\left\vert
\Omega\right\vert ^{2}}\int\limits_{\Omega}\left(  w-u\right)  ^{2}%
\leqslant\frac{1}{\left\vert \Omega\right\vert }C_{18}^{2}\int\limits_{0}%
^{t}\int\limits_{\Omega}\left(  1-z^{-1}\right)  ^{2}\leqslant\nonumber\\
&  \leqslant\frac{1}{\left\vert \Omega\right\vert }C_{12}C_{18}^{2}%
C_{53}C_{54}\max_{s\in(0,\left\Vert v_{0}\right\Vert _{L^{\infty}\left(
\Omega\right)  })}\chi^{2}(s) \label{d3}%
\end{align}
which implies (\ref{wu_bar}). From (\ref{d1}),(\ref{d2}) and (\ref{d3}) the
estimate (\ref{asymp4})\ follows.

We consider now the case $\mu>0$. From (\ref{z}) we obtain%
\begin{align}
&  \frac{1}{2}\frac{d}{dt}\int\limits_{\Omega}\left[  z^{-1}(w-1)^{2}\right]
=-\int\limits_{\Omega}z^{-1}\left\vert \nabla w\right\vert ^{2}-\mu
\int\limits_{\Omega}wz^{-1}(w-1)^{2}+\mu\int\limits_{\Omega}w^{2}%
z^{-1}(1-z^{-1})(w-1)-\nonumber\\
&  -\mu\int\limits_{\Omega}wz^{-1}v(w-1)-\int\limits_{\Omega}\chi
(v)z^{-1}w^{2}v_{t}\leqslant\nonumber\\
&  \leqslant-\int\limits_{\Omega}z^{-1}\left\vert \nabla w\right\vert
^{2}-\frac{\mu}{2}\int\limits_{\Omega}wz^{-1}(w-1)^{2}+\int\limits_{\Omega
}w^{3}z^{-1}(1-z^{-1})^{2}+\int\limits_{\Omega}wz^{-1}v^{2}-C_{57}%
\int\limits_{\Omega}w^{2}v_{t}, \label{sm}%
\end{align}
where
\[
C_{57}=e^{\frac{L_{\chi}}{2}\Vert v_{0}\Vert_{L^{\infty}(\Omega)}^{2}%
+\chi(0)\Vert v_{0}\Vert_{L^{\infty}(\Omega)}}\left(  L_{\chi}\Vert v_{0}%
\Vert_{L^{\infty}(\Omega)}+\chi(0)\right)  .
\]
From the above it follows that%
\begin{align*}
&  \frac{1}{2}\frac{d}{dt}\int\limits_{\Omega}\left[  z^{-1}(w-1)^{2}\right]
+\frac{\mu}{2}\min\left\{  1,a\right\}  \int\limits_{\Omega}z^{-1}%
(w-1)^{2}+\int\limits_{\Omega}z^{-1}\left\vert \nabla w\right\vert
^{2}\leqslant\\
&  \leqslant\int\limits_{\Omega}w^{3}z^{-1}(1-z^{-1})^{2}+\int\limits_{\Omega
}wz^{-1}v^{2}-C_{57}\int\limits_{\Omega}w^{2}v_{t}.
\end{align*}
Integrating the last inequality on $(0,t)$ we have%
\begin{align}
&  \frac{1}{2}\int\limits_{\Omega}\left[  z^{-1}(w-1)^{2}\right]  +\frac{\mu
}{2}\min\left\{  1,a\right\}  \int\limits_{0}^{t}\int\limits_{\Omega}%
z^{-1}(w-1)^{2}+\int\limits_{0}^{t}\int\limits_{\Omega}z^{-1}\left\vert \nabla
w\right\vert ^{2}\leqslant\nonumber\\
&  \leqslant C_{12}C_{18}^{3}\int\limits_{0}^{t}\int\limits_{\Omega}%
(1-z^{-1})^{2}+C_{12}C_{18}\int\limits_{0}^{t}\int\limits_{\Omega}v^{2}%
+C_{18}^{2}C_{57}\left\vert \Omega\right\vert \left\Vert v_{0}\right\Vert
_{L^{\infty}\left(  \Omega\right)  }+\frac{1}{2}\int\limits_{\Omega}%
(w_{0}-1)^{2}e^{\int\limits_{0}^{v_{0}}\chi(s)ds}. \label{sg}%
\end{align}
Taking into account (\ref{sh}) and (\ref{gv1}) with $q=1$, we obtain from
(\ref{sg})%
\begin{equation}
\frac{1}{2}\int\limits_{\Omega}\left[  z^{-1}(w-1)^{2}\right]  +\frac{\mu}%
{2}\min\left\{  1,a\right\}  \int\limits_{0}^{t}\int\limits_{\Omega}%
z^{-1}(w-1)^{2}+\int\limits_{0}^{t}\int\limits_{\Omega}z^{-1}\left\vert \nabla
w\right\vert ^{2}\leqslant C_{58}, \label{sr}%
\end{equation}
where
\[
C_{58}=C_{12}C_{18}C_{54}\left(  C_{12}C_{18}\max_{s\in(0,\left\Vert
v_{0}\right\Vert _{L^{\infty}\left(  \Omega\right)  })}\chi^{2}(s)+C_{53}%
^{-1}\right)  +C_{57}C_{18}^{2}\left\vert \Omega\right\vert \left\Vert
v_{0}\right\Vert _{L^{\infty}\left(  \Omega\right)  }+\frac{1}{2}%
\int\limits_{\Omega}(w_{0}-1)^{2}e^{\int\limits_{0}^{v_{0}}\chi(s)ds}.
\]
The inequality (\ref{sr}) implies%
\begin{align}
\int\limits_{0}^{t}\int\limits_{\Omega}z^{-1}(w-1)^{2}  &  \leqslant
C,\label{ab}\\
\int\limits_{0}^{t}\int\limits_{\Omega}z^{-1}\left\vert \nabla w\right\vert
^{2}  &  \leqslant C. \label{ac}%
\end{align}
From (\ref{sm}) we deduce that%
\begin{align*}
\left\vert \frac{d}{dt}\int\limits_{\Omega}\left[  z^{-1}(w-1)^{2}\right]
\right\vert  &  \leqslant2\int\limits_{\Omega}z^{-1}\left\vert \nabla
w\right\vert ^{2}+2\mu\int\limits_{\Omega}wz^{-1}(w-1)^{2}+2\mu\int
\limits_{\Omega}w^{2}z^{-1}\left\vert (1-z^{-1})(w-1)\right\vert +\\
&  +2\mu\int\limits_{\Omega}wz^{-1}v\left\vert w-1\right\vert -2\int
\limits_{\Omega}\chi(v)z^{-1}w^{2}v_{t}%
\end{align*}
and integrating in time we have
\begin{align*}
&  \int\limits_{0}^{t}\left\vert \frac{d}{dt}\int\limits_{\Omega}\left[
z^{-1}(w-1)^{2}\right]  \right\vert \leqslant2\int\limits_{0}^{t}%
\int\limits_{\Omega}z^{-1}\left\vert \nabla w\right\vert ^{2}+2\mu C_{18}%
\int\limits_{0}^{t}\int\limits_{\Omega}z^{-1}(w-1)^{2}+\mu C_{12}C_{18}%
^{2}\int\limits_{0}^{t}\int\limits_{\Omega}(1-z^{-1})^{2}+\\
&  +\mu C_{18}^{2}\int\limits_{0}^{t}\int\limits_{\Omega}z^{-1}(w-1)^{2}+\mu
C_{12}C_{18}^{2}\int\limits_{0}^{t}\int\limits_{\Omega}v^{2}+\mu
\int\limits_{0}^{t}\int\limits_{\Omega}z^{-1}\left(  w-1\right)  ^{2}-\\
&  -2C_{12}C_{18}^{2}\left(  L_{\chi}\Vert v_{0}\Vert_{L^{\infty}(\Omega
)}+\chi(0)\right)  \int\limits_{0}^{t}\int\limits_{\Omega}v_{t}\leqslant\\
&  \leqslant2\int\limits_{0}^{t}\int\limits_{\Omega}z^{-1}\left\vert \nabla
w\right\vert ^{2}+\mu\left(  C_{18}+1\right)  ^{2}\int\limits_{0}^{t}%
\int\limits_{\Omega}z^{-1}(w-1)^{2}+\mu C_{12}C_{18}^{2}\int\limits_{0}%
^{t}\int\limits_{\Omega}(1-z^{-1})^{2}+\\
&  +\mu C_{12}C_{18}^{2}\int\limits_{0}^{t}\int\limits_{\Omega}v^{2}%
+2C_{12}C_{18}^{2}\left\vert \Omega\right\vert \left(  L_{\chi}\Vert
v_{0}\Vert_{L^{\infty}(\Omega)}+\chi(0)\right)  \Vert v_{0}\Vert_{L^{\infty
}(\Omega)}.
\end{align*}
Using (\ref{gv1}), (\ref{sh}), (\ref{ab}) and (\ref{ac}) we conclude%
\begin{equation}
\int\limits_{0}^{t}\left\vert \frac{d}{dt}\int\limits_{\Omega}\left[
z^{-1}(w-1)^{2}\right]  \right\vert \leqslant C. \label{zw2}%
\end{equation}
The last inequality and (\ref{ab}) imply%
\begin{equation}
\lim_{t\rightarrow\infty}\int\limits_{\Omega}\left[  z^{-1}(w-1)^{2}\right]
=0. \label{zw3}%
\end{equation}
From (\ref{zw3}) and taking into account the following estimate%
\[
\Vert u(\cdot,t)-1\Vert_{L^{2}(\Omega)}^{2}\leqslant e^{\max\limits_{s\in
(0,\left\Vert v_{0}\right\Vert _{L^{\infty}\left(  \Omega\right)  })}\chi
(s)}\int\limits_{\Omega}z^{-1}\left(  w-1\right)  ^{2}+\Vert z^{-1}%
-1\Vert_{L^{2}(\Omega)},
\]
we conclude the proof.
\end{proof}

\section{Acknowledgments}

This work was partially supported by the RTN \textquotedblright Modeling,
Mathematical Methods and Computer Simulation of Tumour Growth and
Therapy\textquotedblright\ (MRTN-CT-2004-503661). The first author was also
partially supported by project DGES (Spain) Grant MTM2007-61755. The authors
would like to thank the Centre de Recerca Matem\`{a}tica, where this paper was
finished, for the invitation to participate to the research programme on
"Mathematical Biology: Modelling and Differential Equations" during February
2009 and for the excellent working conditions provided.

\end{document}